\documentclass[10pt]{amsart}

\usepackage{fullpage}
\usepackage{amsmath, amsfonts, amssymb, amsthm, mathdots, mathrsfs, dsfont}
\usepackage{float}
\usepackage{tikz,tikz-cd}
\usepackage[colorlinks, urlcolor=black]{hyperref}

\hypersetup{citecolor=blue, linkcolor=blue}

\newcommand{\Q}{\mathbb{Q}}

\newcommand{\N}{\mathbb{N}}

\newcommand{\F}{\mathbb{F}}
\newcommand{\A}{\mathbb{A}}
\newcommand{\C}{\mathbb{C}}
\DeclareMathOperator{\sym}{Sym}
\DeclareMathOperator{\Orb}{Orb}
\DeclareMathOperator{\Gal}{Gal}

\DeclareMathOperator{\PSO}{PSO}

\DeclareMathOperator{\PO}{PO}

\DeclareMathOperator{\PGO}{PGO}
\DeclareMathOperator{\POM}{P\Omega}
\DeclareMathOperator{\GO}{GO}
\DeclareMathOperator{\GL}{GL}
\DeclareMathOperator{\PSL}{PSL}
\DeclareMathOperator{\PSU}{PSU}
\DeclareMathOperator{\PGL}{PGL}
\DeclareMathOperator{\PGU}{PGU}
\DeclareMathOperator{\PSp}{PSp}
\DeclareMathOperator{\PGSp}{PGSp}

\newtheorem{theorem}{Theorem}[section]
\newtheorem{definition}[theorem]{Definition}
\newtheorem{proposition}[theorem]{Proposition}

\newtheorem{corollary}[theorem]{Corollary}
\newtheorem{lemma}[theorem]{Lemma}
\newtheorem{remark}[theorem]{Remark}

\newtheorem{example}[theorem]{Example}

\numberwithin{equation}{section}

\author{Luis Gutiérrez Frez
and Adri\'an Zenteno}
\begin{document}

\title{Classifying representations of finite classical groups of Lie type of dimension up to $\ell^4$}

\begin{abstract}
Let $G$ be a finite classical group of Lie type of rank $\ell$, defined over a field of characteristic $p>2$. In this work, we classify the irreducible representations of $G$ whose dimensions are bounded by a constant proportional to $\ell$, and splits into two cases according to $G$ is of type $A_{\ell}$ or not. Furthermore, we  discuss explicit formulas for computing the dimensions of such representations.
 The motivation for this work arises, in part, from a desire to obtain new results on two classical problems concerning Galois representations: the large image conjecture for automorphic Galois representations and the inverse Galois problem. We conclude the paper by giving some remarks on potential implications in these addresses.

\end{abstract}

\maketitle

\section{Introduction}

In this paper we are  particularly interested in the representation theory of the finite classical groups $G$ of Lie type of rank $\ell$.  It is well known that every irreducible representation $V$ of a simple connected, simple algebraic group over a field $K$ is in correspondence with a suitable dominant weight   $\lambda$, called highest weight, of a root system associated to maximal torus $T$ of $G$. In particular, for any dominant weight $\lambda$, there exists an irreducible finite-dimensional  $G(K)$-module $L(\lambda)$ of highest weight $\lambda$. This module is not fully understood for an arbitrary weight $\lambda$. For instance, differences in the dimension of $L(\lambda)$ are expected according to the characteristic of $K$: in characteristic zero, the dimension of $L(\lambda)$ is given by the Weyl's formula, whereas in positive characteristic $p$, a  closed formula for the dimension of $L(\lambda)$ is still far to be known for an arbitrary weight $\lambda$. In this case, the multiplicities involved in computing $\dim L(\lambda)$ depend on $p$, and determining this dependence is a challenging task. In this regard, progress have been established case-by-case according to $G$ and  very special weights $\lambda$ (see \cite{Bur,BW,BGT,Ca,GS,Li,L,Sei}-for instance). With this in mind, we  aim to classify irreducible modules whose dimension is bounded by a given constant.  More precisely, we will say that $\lambda$ is an admissible weight if 
\begin{equation}\label{admis}
 \dim L(\lambda)\leq 
 \begin{cases}
\ell^n/2^n& \text{ if }G \text{ of type  }A_{\ell}\\
\ell^n, & \text{ otherwise}.\\
\end{cases}   
\end{equation}

A complete classification of all admissible weights for arbitrary $n$ turns out to be a major challenge, which makes a case-by-case approach both natural and plausible.
This way was already started by Liebeck in \cite{Li}, where all the admissible weights for $n=2$ were determined. Later, L\"ubeck \cite{L} extended those results in the case of $n=3$, providing, in addition, a series of tables listing representations whose dimension is less than or equal to a bound $M$, depending on the type.  An extension of  such a classification for bounds $\binom{\ell+1}{4}$ or $16\binom{\ell}{4}$ according to the type is given by Mart\'inez  \cite{Mar1}.  

The main result of this paper split into three theorems, in which we classify all admissible 
weights when $n=4$ in Eq \ref{admis}. More precisely, in \S\ref{S.3}, we state our first result -Theorem \ref{PrinA}- concerning the classification of all the admissible weights  when $G$ is of type $A_{\ell}$.  The proof  begins with an arbitrary dominant  weight $\lambda$ and  consists on determining all  the possibilities of its coefficients via a series of reductions until to left reduced to a full list of admissible weights. 
Later, in  \S \ref{S.5}  we  provide the classification of admissible weights for groups $G$ of type $B_{\ell}$, $C_{\ell}$ and $D_{\ell}$. Here, our method implies to obtain a series of bounds for the multiplicities for the weights involved in the classification. In this vein, for  $B_{\ell}$, we have to prove a series of lemmas  on required multiplicities in order to obtain a list of admissible weights -achieved in Proposition \ref{DB}. We then deduce, in Theorem \ref{PrinB},  that  such list exhausts all the admissible weights by means of successive reductions. The section  concludes with Theorem \ref{PrinCD}, where all admissible weights  for the remaining two types are classified by appealing to Appendix -Tables 1 and 2- for the required multiplicities.

As noted above, determining the exact dimension of $L(\lambda)$ in positive characteristic  is significantly more complicated than in characteristic zero, since a closed formula for this dimension is unknown. Here, computing the multiplicity $m_{\lambda}(\mu)$, for any dominant weight $\mu\preceq \lambda$,  requires constructing a $d\times d$-matrix $C^{(\lambda,\mu)}$ and then the multiplicity of $\mu$ is the rank of  the matrix obtained from $C^{(\lambda,\mu)}$ by reduction modulo $p$, where $d$ depends on the Kostant partition function. In particular, determining the precise dependence of the multiplicities involved in $\dim L(\lambda)$  upon reduction modulo $p$ relies heavily on the nature of $\lambda$, and demands substantial effort.
In $\S 5$, we compile and supplement the dimension formulas of $L(\lambda)$ for all our admissible weights $\lambda$ when $p$ is sufficiently large.

Now, let $G_\Q:=\Gal(\overline{\Q}/\Q)$ be the absolute Galois group of $\Q$, and let $\pi$ be a cuspidal, algebraic regular, automorphic representations of $\GL_n(\A_\Q)$, where $\A_\Q$ denotes the ring of adeles of $\Q$. According to global Langlands correspondence \cite{BG}, for each prime $p$ there should exists a $p$-adic Galois representation $\rho_{\pi,p} : G_\Q \rightarrow \GL_n(\overline{\Q}_p)$ associated to $\pi$. A folklore conjecture predicts that the image of the projectivization $\overline{\rho}^{proj}_{\pi,p} : G_\Q \rightarrow \PGL_n(\overline{\F}_p)$ of the residual representation $\overline{\rho}_{\pi,p} : G_\Q \rightarrow \GL_n(\overline{\F}_p)$ of $\rho_{\pi,p}$ should be as large as possible for almost all primes $p$ ({\it i.e.}, all but finitely many), except when  an automorphic reason impedes it.   This statement is often referred to as the (residual) large image conjecture and it has interesting consequences on the Inverse Galois problem for finite groups of Lie type \cite{ADW1, ADW2, Rib, Ze20, Ze24, Ze23}.  It is worth noting that in the few known cases ($n \leq 6$) of this conjecture for $\pi$ polarizable, its proof is done by considering all possible images of $\overline{\rho}^{proj}_{\pi,p}$ given by the classification of the maximal subgroups of $\PSp_n(\F_{p^s})$ or $\POM^\pm_n(\F_{p^s})$. 

In the eighties, Aschbacher \cite{As} classified the maximal subgroups of almost all of the finite almost simple classical groups,   dividing them into two classes.
One of  these consists of absolutely irreducible subgroups that arise (at least for sufficiently large $p$) as representations of finite groups of Lie type. However, this classification is not fine enough to contribute directly to the large image conjecture. The search for a more refined classification is precisely one of the motivations behind this work.
We conclude this paper by showing a weak version of the large image conjecture, and  obtain new weak results on the inverse Galois problem as illustrated in Theorem \ref{princi} and Corollary \ref{Galois}.


\section{Preliminaries on representations of finite simple groups of Lie type}\label{S.2}

In this section, we recall some results about algebraic groups and their defining characteristic representations that will be needed in this paper. 


\subsection{Root Systems} Let $G$ be a connected reductive simple algebraic group over an algebraically closed field $K$ of simply connected type and rank $\ell$.
Let $T$ be a maximal torus of $G$, $X \cong \mathbb{Z}^{\ell}$ its character group and $Y \cong \mathbb{Z}^{\ell}$ its co-character group. Let $\{ \alpha_1, \ldots, \alpha_{\ell} \} \subset X$ be a set of simple roots for $G$ with respect to $T$ and $\Delta$ its associated root system. Let $\alpha_i^\vee \in Y$ be the coroot corresponding to $\alpha_i$ for each $1 \leq i\leq \ell$. 
Viewing $X \otimes \mathbb{R}$ as a Euclidean space we can define the \emph{fundamental weights} $\lambda_1 \ldots, \lambda_{\ell} \in X \otimes \mathbb{R}$ as the dual basis of $\alpha_1^\vee, \ldots, \alpha_{\ell}^\vee$. 

For our purposes, and in accordance with \cite{L}, we shall use the opposite order from that given in  \cite{Bou02,Hum75}. Thus, the Dynkin diagrams for classical groups are given  as follows:

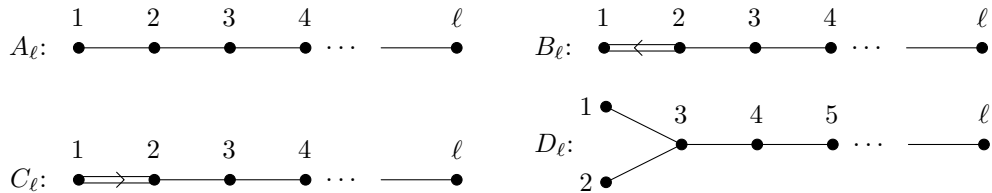
\begin{align*}
\begin{tikzpicture}[line cap=round,line join=round,x=1.0cm,y=1.0cm]
 \draw (0.7,0) node[anchor=east]  {$A_{\ell}$:};
\draw (1,0)-- (4,0);
\draw (5,0)-- (6,0);
\draw (1.2,0.4) node[anchor=east]  {$1$};  
\draw (2.2,0.4) node[anchor=east]  {$2$};  
\draw (3.2,0.4) node[anchor=east]  {$3$};  
\draw (4.2,0.4) node[anchor=east]  {$4$};  
\draw (6.2,0.4) node[anchor=east]  {$\ell$};  
\draw[fill=black] (1,0) circle (2pt);
\draw[fill=black] (2,0) circle (2pt);
\draw[fill=black] (3,0) circle (2pt);
\draw[fill=black] (4,0) circle (2pt);
\draw[color=black] (4.5,0) node {$\dots$};
\draw[fill=black] (6,0) circle (2pt);
\end{tikzpicture}
\quad&\quad
\begin{tikzpicture}[line cap=round,line join=round,x=1.0cm,y=1.0cm]
\draw (0.7,0) node[anchor=east]  {$B_{\ell}$:};
\draw (1,0)-- (4,0);
\draw (5,0)-- (6,0);
\draw[double distance=2pt] (1,0)-- (2,0);
\draw (1.2,0.4) node[anchor=east]  {$1$};  
\draw (2.2,0.4) node[anchor=east]  {$2$};  
\draw (3.2,0.4) node[anchor=east]  {$3$};  
\draw (4.2,0.4) node[anchor=east]  {$4$};  
\draw (6.2,0.4) node[anchor=east]  {$\ell$};  
\draw[fill=black] (1,0) circle (2pt);
\draw[fill=black] (2,0) circle (2pt);
\draw[fill=black] (3,0) circle (2pt);
\draw[fill=black] (4,0) circle (2pt);
\draw[color=black] (4.5,0) node {$\dots$};
\draw[fill=black] (6,0) circle (2pt);
\draw (1.5,-0.1)-- (1.4,0)-- (1.5,0.1);
\end{tikzpicture}\\
\begin{tikzpicture}[line cap=round,line join=round,x=1.0cm,y=1.0cm]
\draw (0.7,0) node[anchor=east]  {$C_{\ell}$:};
\draw (1,0)-- (4,0);
\draw (5,0)-- (6,0);
\draw[double distance=2pt] (1,0)-- (2,0);
\draw (1.5,0.1)-- (1.6,0)-- (1.5,-0.1);
\draw (1.2,0.4) node[anchor=east]  {$1$};  
\draw (2.2,0.4) node[anchor=east]  {$2$};  
\draw (3.2,0.4) node[anchor=east]  {$3$};  
\draw (4.2,0.4) node[anchor=east]  {$4$};  
\draw (6.2,0.4) node[anchor=east]  {$\ell$};  
\draw[fill=black] (1,0) circle (2pt);
\draw[fill=black] (2,0) circle (2pt);
\draw[fill=black] (3,0) circle (2pt);
\draw[fill=black] (4,0) circle (2pt);
\draw[color=black] (4.5,0) node {$\dots$};
\draw[fill=black] (6,0) circle (2pt);
\end{tikzpicture}
\quad&\quad
\begin{tikzpicture}[line cap=round,line join=round,x=1.0cm,y=1.0cm]
\draw (0.7,0) node[anchor=east]  {$D_{\ell}$: };
\draw (2,0)-- (4,0);
\draw (5,0)-- (6,0);
\draw (1,0.5)-- (2,0)--(1,-0.5);
\draw (0.95,0.5) node[anchor=east]  {$1$};  
\draw (0.95,-0.5) node[anchor=east]  {$2$};  
\draw (2.2,0.4) node[anchor=east]  {$3$};  
\draw (3.2,0.4) node[anchor=east]  {$4$};  
\draw (4.2,0.4) node[anchor=east]  {$5$};  
\draw (6.2,0.4) node[anchor=east]  {$\ell$};  
\draw[fill=black] (1,0.5) circle (2pt);
\draw[fill=black] (1,-0.5) circle (2pt);
\draw[fill=black] (2,0) circle (2pt);
\draw[fill=black] (3,0) circle (2pt);
\draw[fill=black] (4,0) circle (2pt);
\draw[color=black] (4.5,0) node {$\dots$};
\draw[fill=black] (6,0) circle (2pt);
\end{tikzpicture}
\end{align*}
We also recall that  each $\alpha\in\Delta$ determines a reflection $\sigma_{\alpha}$ respect to the hyperplane $P_{\alpha}=\{\beta\in E:(\beta,\alpha)=0  \}$, 
where $(E,(-,-))$ is the underlying  inner product real vector space of finite dimensional associated to $\Delta$.
 The Weyl group $\mathcal{W}$ of $\Delta$ is the subgroup  of the orthogonal group on $E$ generated by the reflections $\sigma_{\alpha_i}$.
Thus, one observes that this group acts on $\Delta$ by permutation.  In fact,  $\mathcal{W}$ is generated by the reflections 
$\sigma_{\alpha_i}$ with $1\leq i\leq \ell$ by appealing to \cite[Prop. 2.62]{Kna}.

Furthermore, one obtains a partial ordering in $X$ by declaring that a weight $\lambda$ is smaller than a weight $\lambda'$ (denoted by $\lambda \prec \lambda'$) if and only if $\lambda'-\lambda$ is a non-negative linear combination of simple roots. We shall also use $\preceq$ in an obvious way.  In addition,   a weight $\lambda \in X$ is dominant if it is a non-negative linear combination of the fundamental weights. The group $\mathcal{W}$ acts on the set of the dominant weights by means of its action on the simple roots $\alpha_1,\cdots,\alpha_{\ell}$. Given a dominant weight $\lambda$,  we denote by 
 $\Orb_\mathcal{W}(\lambda)$ and $\rm{Stab}_{\mathcal{W}}(\lambda)$  the orbit and the stabilizer of $\lambda$ under the $\mathcal{W}$-action, respectively. \\

It is worth recalling the following.

\begin{remark} \label{stab} 
The stabilizer subgroup $\rm{Stab}_{\mathcal{W}}(\lambda)$ of a dominant  weight $\lambda=a_1\lambda_1+a_2\lambda_2+\cdots+a_{\ell}\lambda_{\ell}$  is the parabolic subgroup of $G$ generated by the reflections along the simple roots $\alpha_i$ such that $a_i=0$.    
\end{remark}
It  will be used  implicitly many times along the next sections mentioning it only the  first  instances we use it.

\subsection{Representations} 
Henceforth, we assume that $K = \overline{\F}_p$ and let $V$ be a finite-dimensional $G(K)$-module. A classical result tells us 
\begin{equation}
V=\bigoplus_{\mu\in X} V_{\mu},    
\end{equation}
where $V_\mu := \{ v \in V: vt = \mu(t)v \mbox{ for all }t \in T \}$, called weight space, for any dominant weight $\mu\in X$.
Each  $\mu \in X$ which satisfies $V_\mu \neq \{ 0 \}$ is called a weight of $V$. From a result of Chevalley (see Section 31.3 of \cite{Hum75}) we have that, when $V$ is irreducible, there exists a unique $T$-stable $1$-dimensional subspace spanned by a unique maximal vector of dominant weight $\lambda$, is called the \emph{highest weight} of $V$.
Moreover, we have that each irreducible $G(K)$-module $V$ is determined (up to isomorphism) by its highest weight and for each dominant weight $\lambda \in X$  there exists an irreducible $G(K)$-module $L(\lambda)$ with highest weight $\lambda$.
A way to construct $L(\lambda)$ is to start with a simple module 
$V(\lambda)_{\mathbb{C}}$ over $\C$ for the corresponding semisimple Lie algebra $\mathfrak{g}$. The theory in characteristic zero provides the existence of a special $\mathbb{Z}$-form $V(\lambda)_{\mathbb{Z}}\subset V(\lambda)_{\mathbb{C}}$, which allows reduction modulo $p$, such that $G$ acts naturally on the so-called Weyl module $V(\lambda):=K \otimes V(\lambda)_{\mathbb{Z}}$. Thus, the irreducible $G$-module 
 $L(\lambda)$ is given by
\[
L(\lambda):=V(\lambda)/\rm{rad}(\lambda),
\]
where $\rm{rad}(\lambda)$ is the unique irreducible maximal submodule of $V(\lambda)$.

Let $F$ be the canonical Frobenius automorphism of $K$ given by $x \mapsto x^p$. Then, given a $G(K)$-module $L$ and $i\geq 0$, we denote by $L^{(i)}$ the $G(K)$-module obtained by composing the $G$-action on $L$ with the $i$-th power of $F$. On the other hand, we say that a dominant weight $\lambda = a_1 \lambda_1 + \cdots + a_\ell \lambda_\ell \in X$ is $p$-\emph{restricted} if for every $i=1, \ldots, \ell$ we have $0 \leq a_i \leq p-1$, and we may introduce the Twisted Tensor Product Theorem of Steinberg \cite{Ste63}.

\begin{theorem}\label{tutu} Let $\lambda \in X$ be a dominant weight. Then, there exist $r \in \N$ and $p$-restricted dominant weights $\lambda_0, \cdots, \lambda_r \in X$ such that as $G(K)$-modules,
\[
L(\lambda) = L(\lambda_0 + p \lambda_1 + \cdots +p^r \lambda_r) \cong L(\lambda_0) \otimes L(\lambda_1)^{(1)} \otimes \cdots \otimes L(\lambda_r)^{(r)}.
\]
\end{theorem} 

Then, the study of all irreducible $G(K)$-modules is reduced to that of the finitely many ones having $p$-restricted highest weights. In this vein,  computing the dimension $L(\lambda)$ has been a theme of profound importance and multiple articles have  contributed to obtain  a better understanding of that. Next, we give a brief account of this matter.

\subsection{Dimension and Multiplicity}
We shall continue by presenting the well-known Freudenthal formula in characteristic zero. To this end, let $V_{\mathbb{C}}$ be a finite-dimensional irreducible $G$-module with highest weight $\lambda$. 
We need to introduce a special  weight which  appears in the computation of multiplicities as we will see here below. Specifically,   the  strongly dominant weight $\delta$ is defined by 
    \[
    \delta:=\frac{1}{2}\sum_{\alpha\succ 0}\alpha=\sum_{i=1}^{\ell}\lambda_i,
    \]
    following \cite[\S 10.2 and \S13.3]{Hum}.

 We define the multiplicity of $\mu$ in $V$, for any dominant weight $\mu \preceq \lambda$, as
\[
m_{\lambda}(\mu) = \dim (V_{\mathbb{C}})_\mu  
\]
which can be calculated by using the recursive   formula 
\cite[\S 22.3]{Hum}:
\begin{equation}\label{FreudForm}
m_\lambda (\mu) = \frac{2 \sum_{\alpha \in \Phi_+}  \sum_{i=1}^{\infty} (\mu + i \alpha, \alpha) m_\lambda (\mu + i \alpha)}{ \parallel \lambda + \delta \parallel^2 + \parallel \mu + \delta \parallel^2}
= \frac{2 \sum_{\alpha \in \Phi_+}  \sum_{i=1}^{\infty} (\mu + i \alpha, \alpha) m_\lambda (\mu + i \alpha)}{(\lambda,\lambda)-(\mu,\mu)+2(\lambda-\mu,\delta)},
\end{equation}
where $\Phi_+$ denotes the set of positive roots of $\Delta$ and $m_\lambda(\lambda) = 1$.

It follows from \cite[(1.10)]{Sei} (see also \cite[Proposition 2.4.3]{CTesis}) that the dimension of $V_{\mathbb{C}}$ is

\begin{equation}\label{dimf}
    \dim V_{\mathbb{C}} = \sum _{\mu \in \Pi(\lambda)} m_{\lambda}(\mu) l(\mu),
\end{equation}
where $\Pi(\lambda) := \{ \mu : \mu \mbox{ is a dominant weight such that }\mu \preceq \lambda\}$ and $l(\mu)$ denotes the length of the $\mathcal{W}$-orbit $\Orb_\mathcal{W}(\mu)$. 

Nevertheless, as we mentioned in the Introduction, determining the dimension of $L(\lambda)$ in positive characteristic is much more complicated. In particular, knowing  multiplicity $m_{\lambda}(\mu)$  is to compute the rank of  the matrix obtained from $C^{(\lambda,\mu)}$ by reduction modulo $p$. Now, given that   those calculations are  very lengthy even in low rank, computational programs have been a useful tool  since the 1970s. In addition, we may observe that the multiplicities obtained by reduction modulo $p$  are not, in general, the same, the appearance of the weights $\mu$ also demands attention. For instance, Premet  established that the set of weights for the modules $V$, $L(\lambda)$ is the same if $G$ is of type $B_{\ell}$ and $C_{\ell}$ when $p>3$ and $\lambda$ is $p$-restricted,  which  are the central ones for us as we will see below. 

Thereby, it follows that 
\begin{equation}\label{dimf1}
    \dim L(\lambda) \leq \sum _{\mu \in \Pi(\lambda)} m_{\lambda}(\mu) l(\mu) = \dim V_{\mathbb{C}}.
\end{equation}

We finish this section by specifying the notion of {\it admissible weight} mentioned in the Introduction. 

\begin{definition}\label{ADM} Let $G$ be a connected reductive simple algebraic group of rank $\ell$ of type $A,B,C,D$ over $\overline{\F}_p$ of simply connected type and $\lambda$ be a dominant weight.
We say that  $\lambda$ is admissible if 
\[
\dim L(\lambda)\leq 
 \begin{cases}
\ell^4/16, & \text{ if }G \text{ is of type  }A_{\ell}\\
\ell^4, & \text{ otherwise}.\\
\end{cases}
\]
\end{definition}
The next three sections will be devoted to determining all the admissible weights  for all finite classical groups in  their defining characteristic.


\section{Representations of finite classical groups of  type $A_\ell$}\label{S.3}

In this section we will determine all the admissible weights of type $A_\ell$ together with the exact dimension of their associated irreducible representations.

Let $G$ be a finite classical group of type $A_{\ell}$. We shall begin by fixing a set   of simple roots  of type $A_{\ell}$ and its fundamental weights associated. 
To this end, let $\{\alpha_1,\alpha_2,\cdots,\alpha_{\ell}\}$ be the simple root system of type $A_\ell$ given by
$$
\alpha_i:=e_i-e_{i+1},\quad\quad 1\leq i\leq \ell.
$$
We can normalize as follows: $(\alpha_i, \alpha_i)= 1$, $(\alpha_i, \alpha_{i+1})=-1/2$ and $(\alpha_i, \alpha_{i+k})=0$, $k >1$. 
Then, the  fundamental  weights are given by 
\[
\lambda_i=\frac{1}{\ell+1}[(\ell-i+1)\alpha_1+2(\ell-i+1)\alpha_2+i(\ell-i+1)\alpha_i+i(\ell-i)\alpha_{i+1}+
\cdots+i\alpha_{\ell}],
\]
for $1 \leq i\leq \ell$.

We commence by obtaining the following admissible weights.

\begin{lemma}\label{4AS}

Let $G$ be as above with $\ell \geq 15$. Then, the weights 
\[
3\lambda,\; 4\lambda,\; \lambda_3, \;\lambda_4,
\;\lambda_1+\lambda_2,\;\lambda_1+\lambda_{\ell-1},\; \lambda_1+2\lambda_{\ell}
\]

(and its duals) are admissible.
\end{lemma}
\begin{proof}
The weights $3\lambda_{\ell}$ and  $4\lambda_{\ell}$  are admissible, since  
\[
\dim L(a\lambda_{\ell})=\binom{\ell+a}{a}<\frac{\ell^4}{4},
\]
for $a=3,4$ and  any $p>3$, by appeal to \cite[Lemma 2.3]{Ca}. The admissibility of $ \lambda_3, \;\lambda_4$ can be deduced by using  item 2 of \cite[Lemma 2.3]{Ca} independently of $p$. On the other hand, item 3, for $2\leq j\leq \ell$ in loc. cit., tells us that 
\begin{equation}\label{noadm}
\dim L(\lambda_1+\lambda_{j})=j\binom{\ell+2}{j+1}-
\epsilon_{p}(j+1)\binom{\ell+1}{j+1},
\end{equation}

where $\epsilon_{p}(j+1)$ is $1$ if $p$ divides $j+1$ and $0$ otherwise. An evident  lower bound for the dimension is when $p$ does not divide $j+1$, from which we infer the admissibility of $\lambda_1+\lambda_{j}$ for $j=2$
and $j=\ell-1$. Finally, $\lambda_1+2\lambda_{\ell}$ is also admissible by appealing to $V(\lambda_1+2\lambda_{\ell})$ is contained in 
$V(\lambda_1)\otimes V(2\lambda_\ell)$ whose dimension is clearly less that $\ell^4/16$.
\end{proof}

\begin{remark}\label{pesofaltante}
An explicit formula for the dimension of $L(\lambda_1 + 2 \lambda_\ell)$ 
    \[
    \dim L(\lambda_1 + 2 \lambda_\ell) = \dim V(\lambda_1 + 2 \lambda_\ell) - \epsilon_p(j+2) L(\lambda_1) = (\ell+1) \left( \frac{\ell(\ell+3)}{2} - \epsilon_p(j+2) \right)
    \]
can be found in \cite[Lemma 2.4.9]{CTesis}.
\end{remark}

In contrast with the example above, we observe.
\begin{example}\label{B21}
For all $\ell \geq 12$, the weight $\lambda=\lambda_{\ell-2}+\lambda_{\ell}$ is not admissible for any $p$. Concretely, 
by considering $j=3$ in  Eq. (\ref{noadm}). we may infer that 
\begin{equation*}
\dim L(\lambda_1+\lambda_{j})=3\binom{\ell+2}{4}-
\binom{\ell+1}{4}=\frac{1}{12}(\ell-1)\ell(\ell+1)(\ell+4)>\frac{\ell^4}{16}.
\end{equation*}

\end{example}

Now we prove the main result of this section.
\begin{theorem}\label{PrinA}
Let  $G$ be a finite group of Lie type $A_\ell$ over $\mathbb{\overline{F}}_p$. If $\ell \geq 15$ then 
the following list exhausts the admissible $p$-restricted weights.

\begin{center}
\begin{tabular}{ l r r}
$\lambda$ & dimension of $L(\lambda)$ & characteristic \\\hline
 $0$ & $1$ & all \\
 $\lambda_1$, $\lambda_\ell$ & $\ell+1$ & all \\
$\lambda_2$, $\lambda_{\ell-1}$ & $\ell(\ell+1)/2$ & all \\
$\lambda_3$, $\lambda_{\ell-2}$ & $(\ell-1)\ell(\ell+1)/6$ & all \\
$\lambda_4$, $\lambda_{\ell-3}$  & $(\ell-2)(\ell-1)\ell(\ell+1)/24$ & all \\
$2 \lambda_1$, $ 2\lambda_{\ell}$ & $(\ell+2)(\ell+1)/2$ & $p \neq 2$ \\
$3 \lambda_1$, $ 3\lambda_{\ell}$ & $(\ell+3)(\ell+2)(\ell+1)/6$ & $p \neq 2,3$ \\
$\lambda_1 + \lambda_\ell$ & $\ell (\ell+2)$ & $p \nmid \ell+1$ \\
$\lambda_1 + \lambda_\ell$ &  $\ell (\ell+2)-1$ &  $p \mid \ell+1$ \\
$\lambda_1 + \lambda_{\ell-1}$, $\lambda_2 + \lambda_{\ell}$  & $(\ell+2)(\ell+1)(\ell-1)/2$ & $p \nmid \ell$ \\
$\lambda_1 + \lambda_{\ell-1}$, $\lambda_2 + \lambda_{\ell}$  &  $(\ell+2)(\ell+1)(\ell-1)/2 - (\ell+1) $  & $p\mid \ell $ \\
$\lambda_1 + 2\lambda_{\ell}$, $2 \lambda_1 + \lambda_{\ell}$ & $(\ell+1)\ell(\ell+3)/2$ & $p \neq 2$, \; $p \nmid \ell+2$ \\
$\lambda_1 + 2\lambda_{\ell}$, $2 \lambda_1 + \lambda_{\ell}$ &  $(\ell+1)\ell(\ell+3)/2-(\ell+1)$ & $p \neq 2$, \; $p \mid \ell+2$ \\
$\lambda_1 + \lambda_2$, $\lambda_{\ell-1} + \lambda_{\ell}$ & $\ell (\ell+1) (\ell+2)/3$ & $p \neq 3$ \\
$\lambda_1 + \lambda_2$, $\lambda_{\ell-1} + \lambda_{\ell}$ & $\ell (\ell+1) (\ell+2)/3- (\ell-1) \ell (\ell+1)/6$ & $p=3$ \\
\hline
\end{tabular}
\end{center}
In addition for $\ell\geq 16$ and $p \neq 2,3$, the weights $4\lambda_1$ and $4\lambda_{\ell}$ are also admissible and $p$-restricted, and their associated irreducible representations have dimension  $(\ell+4)(\ell+3)(\ell+2)(\ell+1) /24$.
\end{theorem}
Notice first that all weights in the table above are admissible by appealing to Proposition \ref{4AS} and \cite[Proposition 3.2]{L}. We therefore need to prove that the  aforementioned list exhausts the admissible weights. To this end,  we shall proceed by a series of steps to  reduce the possibilities according to whether the coefficients are zeroes or not. Explicitly, the proof  begins with.
\begin{proof}
Let $\lambda=a_1\lambda_1+a_2\lambda_2+\cdots+a_{\ell}\lambda_{\ell}$ be a dominant weight. 
The first step, Nullity of Coefficients, consists on identifying all indexes $i$ such that $a_i=0$. Each of the other reductions will be associated to a nonzero coefficient $a_i$, where we determine all the possibilities to the other nonzero coefficients $a_j$ reducing exclusively to admissible weights. Explicitly, we start with. \\
{\it Nullity of coefficients $-$ }  Claim: $a_i \neq 0$ only if $i\leq4$ or $i\geq \ell-3$. \\
Suppose $a_{i}\neq 0$ for some $i$ that satisfy $5\leq i\leq  \ell-4$. This implies that   the $\mathcal{W}$-stabilizer of $\lambda$ is contained in a reflection group of type  $A_{i-1}\times A_{\ell-i}$, by using Remark \ref{stab}. Thus,
\[
L(\lambda)) \geq l(\lambda)=  \frac{(\ell+1)!}{i!(\ell+1-i)!}.
\]

We shall prove that $(\ell+1)!/i!(\ell+1-i)!>\ell^4/16$,  for any $5\leq i\leq  \ell-4$.
To do this, assume firstly that $i\leq \ell+1-i$. In this case, we need to show that 
\[
\frac{(\ell+2-i)(\ell+3-i)\cdots(\ell-3)(\ell-2)(\ell-1)\ell(\ell+1)   }{i!}>\frac{\ell^4}{16}.
\]
We shall first verify the inequality for $i=5$, which is equivalent to  
\begin{align*}
\frac{(\ell-3)(\ell-2)(\ell-1)(\ell+1)   }{3\cdot 5}>\frac{\ell^3}{2}.
\end{align*}

Hence, for $\ell\geq 15$, we can observe  that
\[
\frac{(\ell-3)(\ell-2)(\ell-1)(\ell+1)   }{3\cdot 5}>  (\ell-3)(\ell-2)(\ell-1)
\geq \ell^3  (1-3/15) (1-2/15)(1-1/15)>\frac{\ell^3}{2}
\]
by appealing to $(1-3/\ell) (1-2/\ell)(1-1/\ell)>1/2$, since  $1-j/\ell \geq 1-j/15$ for $1\leq j\leq 3$. Thus, the case $i=5$ follows by multiplying  by $\ell/8$. 
The general  case $i\leq \ell+1-i$ is then  deduced from the base case $i=5$ by noting that 
\[
\frac{(\ell+2-i)(\ell+3-i)\cdots(\ell-4 )(\ell-3)(\ell-2)(\ell-1)\ell(\ell+1)   }{i(i-1)\cdots 6\cdot 5!}
\]
\[
= \left(\frac{(\ell+2-i)(\ell+2-(i+1))\cdots(\ell+2-6 )}{i(i-1)\cdots 6}\right)
\left(\frac{(\ell-3)(\ell-2)(\ell-1)\ell(\ell+1)}{5!}\right)
\]
and  the first parenthesis is greater than 1 by using the condition $\ell\geq 15$. Finally, the case $\ell+1-i <i\leq \ell$ follows by using  a symmetric argument. Therefore, we reduce to
\[
\lambda=a_{1}\lambda_{1}+a_{2}\lambda_{2}+a_{3}\lambda_{3}+a_{4}\lambda_{4}+a_{\ell-3}\lambda_{\ell-3}+a_{\ell-2}\lambda_{\ell-2}+
a_{\ell-1}\lambda_{\ell-1}+a_{\ell}\lambda_{\ell}.
\]


{\it Liberty of  coefficients $-$} Assuming that $a_{\ell-3}\neq 0$ we shall show that $a_{\ell-3}= 1$  and  $ a_i =0$ for any  $i \leq 4$ and $i \geq\ell-2$. In fact, 
if  $a_{\ell-j}\neq 0$  for all $j=0,1,2$  then the $\mathcal{W}$-stabilizer of 
$\lambda$ is contained in a reflection group of type $A_{\ell-4}$. We then see  that 
\begin{equation}\label{A_4}
\dim L(\lambda) \geq l(\lambda)= \frac{(\ell+1)!}{(\ell-3)!}= (\ell-2)(\ell-1)\ell(\ell+1)>\frac{\ell^4}{16}.
\end{equation}

On the other hand, when exactly two coefficients  $a_{\ell-j}$  are nonzero,  the $\mathcal{W}$-stabilizer of 
$\lambda$ is contained in a reflection group of type $A_{\ell-3}\times A_2$. An argument as above 
shows that $\dim L(\lambda) >\ell^4/16$. Next,  if just one coefficient $a_{\ell-j}$ is nonzero then  the $\mathcal{W}$-stabilizer of 
$\lambda$ is contained in a reflection group of type  $A_{\ell-4}\times A_1\times A_1$ or $A_{\ell-4}\times A_2$ according to $j=1$ or not. Thus, we may rule out this possibility of $\lambda$. It follows that  $a_{\ell-j}= 0$  for all $j=0,1,2$ and   by a symmetry argument one checks  that $a_{i}= 0$  for all $i=1,2,3$. Thus, we reduce to  $\lambda = a_{\ell-3} \lambda_{\ell-3}$.
We will now prove that  $a_{\ell-3} =1$ by supposing the contrary {\it i.e.,}  $a_{\ell-3} \geq 2$. Then, we  observe that the $(\ell-4)$-coefficient of 
 \[
 \lambda-\alpha_{\ell-3} = a_{\ell-3} \lambda_{\ell-3} - (-\lambda_{\ell-4}+2 \lambda_{\ell-3} - \lambda_{\ell-2}) = \lambda_{\ell-4}+ (a_{\ell-3}-2) \lambda_{\ell-3} + \lambda_{\ell-2} 
 \]
 is one, which is not possible  by Nullity of Coefficients. Whence $a_{\ell-3} =1$ and then  $\lambda=\lambda_{\ell-3}$ is admissible by Proposition \ref{4AS}.
 

Now assume that $a_{\ell-3}=0$ and $a_{\ell-2}\neq 0$, that is, 
\[
\lambda=a_{1}\lambda_{1}+a_{2}\lambda_{2}+a_{3}\lambda_{3}+a_{4}\lambda_{4}+a_{\ell-2}\lambda_{\ell-2}+
a_{\ell-1}\lambda_{\ell-1}+a_{\ell}\lambda_{\ell}.
\]

We shall prove that  $a_{i} =0$ for $1\leq i\leq 4$.  Suppose first $a_1\neq0$ and $a_2=a_3=a_4=0$.  We note that  the $\mathcal{W}$-stabilizer of $\lambda$ is contained in a reflection group of type $A_{\ell-4} $ or $A_{\ell-4}\times A_1$ or  $A_{\ell-4}\times A_2$ 
according to $a_{\ell-1}\neq 0\neq a_{\ell}$ or not. Then, reasoning as above, we can deduce that  $\dim L(\lambda)>\ell^4/16$.  Similarly each possibility of some  
 $a_{i} \neq 0$ for $1\leq i\leq 4$ discards $\lambda$ as an admissible weight, and thus 
\[\lambda=a_{\ell-2}\lambda_{\ell-2}+
a_{\ell-1}\lambda_{\ell-1}+a_{\ell}\lambda_{\ell}.
\]
We now show that  $a_{\ell-2} \leq 1$. Otherwise, assuming that $a_{\ell-2} \geq 2$, one observes that the coefficients  at $(\ell-3)$ and $(\ell-1)$ of the dominant weight $\lambda-a_{\ell-2}$ are nonzero, but this is not possible by reduction above. Hence, we reduce our analysis  to study 
\[
 \lambda=\lambda_{\ell-2}+a_{\ell-1}\lambda_{\ell-1}+a_{\ell}\lambda_{\ell}\quad\text{or }\quad \lambda=a_{\ell-1}\lambda_{\ell-1}+a_{\ell}\lambda_{\ell}.
 \]
When $a_{\ell-2}=1$ and $a_{\ell-1}\geq 2$ then the term  $2\lambda_{\ell-2}$ appears in  $\lambda-\alpha_{\ell-1}$, which  is impossible by applying the reasoning above. So, any $\lambda$ of the form $\lambda_{\ell-2}+a_{\ell-1}\lambda_{\ell-1}+a_{\ell}\lambda_{\ell}$ with $a_{\ell-1}\geq 2$ is not admissible and  we reduce us to $\lambda$ of the form  $\lambda=\lambda_{\ell-2}+\lambda_{\ell-1}+a_{\ell}\lambda_{\ell}$ or $\lambda=\lambda_{\ell-2}+a_{\ell-1}\lambda_{\ell-1}$.  For the case     $\lambda=\lambda_{\ell-2}+\lambda_{\ell-1}+a_{\ell}\lambda_{\ell}$, if 
  $a_{\ell}\geq 2$ we may  observe similarly that the term $2\lambda_{\ell-1}$
 appears in the dominant weight $\lambda-\alpha_{\ell-1}$ which is a contradiction  by iterating the previous reasoning.  This reduces our analysis to study the dominant weights $\lambda_{\ell-2}+\lambda_{\ell-1}+\lambda_{\ell}$ and $\lambda_{\ell-2}+\lambda_{\ell-1}$. 
In the same vein, we discard the possibility of a dominant weight of the form $\lambda_{\ell-2}+a_{\ell}\lambda_{\ell}$ with $a_{\ell}\geq 4$.
The cases $a_{\ell}=2,3$ are reduced again to analyze the dominant weights $\lambda_{\ell-2}+\lambda_{\ell-1}+\lambda_{\ell}$ and $\lambda_{\ell-2}+\lambda_{\ell-1}$. Finally, when  $\lambda=\lambda_{\ell-2}+a_{\ell}\lambda_{\ell}$ with $a_{\ell}\geq 2$ cannot be admissible by using similar arguments as  above. Therefore, we only need to determine the dimensions of $L(\lambda)$ for the following dominant weights 
\[
\lambda_{\ell-2}+\lambda_{\ell-1},\quad
\lambda_{\ell-2}+\lambda_{\ell-1}+\lambda_{\ell},
 \quad\lambda_{\ell-2}+\lambda_{\ell},  \quad\lambda_{\ell-2}.
\]

In the first case we see that $\lambda_{\ell-3}+\lambda_{\ell} \prec \lambda_{\ell-2}+\lambda_{\ell-1}$. Thus, we can deduce as before that 
\begin{align*}
\dim L(\lambda_{\ell-2}+\lambda_{\ell-1}) \geq 
l(\lambda_{\ell-3}+\lambda_{\ell}) &=\frac{\ell^4(1-2/\ell)(1-1/\ell)(1+1/\ell)}{6}
>\frac{\ell^4(1-2/\ell)(1-1/\ell)}{6}
>\ell^4/12>\ell^4/16,   
\end{align*}
since $(1-2/\ell)(1-1/\ell)\geq ((1-2/15)(1-1/15)>1/2$ for any $\ell\geq 15$. 
 Note that $\dim L(\lambda_{\ell-2}+\lambda_{\ell-1})$ can be computed by using \cite[Lemma 2.3.3]{Ca}, however as can be seen in our argument such a strong result is not necessary for our purposes here.
In the same way, we can rule out the second case by using that $\lambda_{\ell-3}+2\lambda_{\ell} 
\prec \lambda_{\ell-2}+\lambda_{\ell-1}+\lambda_{\ell}$. 
Note that $\lambda_{\ell-2}+\lambda_{\ell}$  is not admissible by appealing to Example \ref{B21} and $\lambda_{\ell-2}$ is admissible by Proposition \ref{4AS}.
\\
Next, suppose $a_{\ell-3}=a_{\ell-2}=0$ and $a_{\ell-1}\neq 0$:
\[
\lambda=a_{1}\lambda_{1}+a_{2}\lambda_{2}+a_{3}\lambda_{3}+a_{4}\lambda_{4}+
a_{\ell-1}\lambda_{\ell-1}+a_{\ell}\lambda_{\ell}.
\]

If $a_{\ell-1}\geq 2 $ we note again that $\lambda_{\ell-2}$ and $\lambda_{\ell}$ appear in the writing of the dominant weight $\lambda-\alpha_{\ell-1}$, which is impossible as was shown previously. Thus, we reduce to  $a_{\ell-1}=1$, {\it i.e.},  
\[
\lambda=a_{1}\lambda_{1}+a_{2}\lambda_{2}+a_{3}\lambda_{3}+a_{4}\lambda_{4}+\lambda_{\ell-1}+a_{\ell}\lambda_{\ell}.
\]
In this case we may prove that $a_i=0$ for
$2\leq i\leq 4$. Here, many  cases arise for the $\mathcal{W}$-stabilizer of $\lambda$ according to the vanish for some of  these $a_i$ ($2\leq i\leq 4)$.
For instance, if $a_2\neq 0$ and 
$a_3=a_4=0$ we see that $\mathcal{W}$-stabilizer of $\lambda$ is contained in a reflection group of type
$ A_1\times A_{\ell-4}$  or $A_{\ell-4}$ according to $a_1=0$ or not. Thus, we may verify that $\dim L(\lambda)>\ell^4/16$, {\it i.e.,} $\lambda$ is not admissible.
Indeed, in any case we may deduce that the dimension $L(\lambda)$ is bigger than the bound $\ell^4/16$. This implies $$\lambda=a_1\lambda_1+\lambda_{\ell-1}+a_{\ell}\lambda_{\ell}.$$

Reasoning as above, one verifies that $ a_1$ and $a_{\ell}$ cannot be  simultaneously nonzero.
In fact, assume first that $ a_1\geq  2$ we obtain the dominant weight 
\[
\lambda-\alpha_1=(a_1-2)\lambda_1+\lambda_2+\lambda_{\ell-1}+a_{\ell}\lambda_{\ell},
\]
which is ruled out by  computations as above. On the other hand,  if $a_{\ell}\geq 2$ then 
 $2\lambda_{\ell-1}$ appears in the writing of $\lambda-\alpha_{\ell}$, which is impossible too. Moreover, we may discard $\lambda=\lambda_1+\lambda_{\ell-1}+\lambda_{\ell}$ by arguing  as in Example \ref{B21}.
This lefts  us  the dominant  weights 
$\lambda_1+\lambda_{\ell-1}$, $\lambda_{\ell-1}+\lambda_{\ell}$ and $\lambda_{\ell-1}$ which are admissible by Proposition \ref{4AS}.  
Finally suppose that $a_{\ell-3}=a_{\ell-2}=a_{\ell-1}=0$ and $a_{\ell}\neq 0$, {\it i.e.}, we will deal with 
\[
\lambda=a_{1}\lambda_{1}+a_{2}\lambda_{2}+a_{3}\lambda_{3}+a_{4}\lambda_{4}+a_{\ell}\lambda_{\ell}.
\]
We may verify, similarly as above, that $a_3=a_4=0$,  reducing our analysis to  $$\lambda=a_1\lambda_1+a_2\lambda_{2}+a_{\ell}\lambda_{\ell}.$$
Observe that  if $a_2 \neq 0$ then $a_{1} = 0$ by using a dual reasoning to the previous case. This implies that $\lambda = \lambda_2+a_{\ell}\lambda_{\ell}$, which gives just one new admissible weight $\lambda=\lambda_2+\lambda_{\ell}$. Now, if $a_2 = 0$ then we have that $\lambda=a_1\lambda_1+a_{\ell}\lambda_{\ell}$ and we get the new admissible weights $\lambda_1+\lambda_{\ell}$, $2\lambda_1+\lambda_{\ell}$, $\lambda_1+2\lambda_{\ell}$ and $a_{\ell}\lambda_{\ell}$ for $a_{\ell}=1,2,3,4$.
This concludes the proof concerning admissibility. Finally, Proposition \ref{4AS}, Remark \ref{pesofaltante}, \cite[Theorem 5.1]{L} and \cite[Lemma 2.3]{Ca} give us the exact dimensions.
\end{proof}


\section{Multiplicities and dimension in finite classical groups}\label{S.5}

In this section we shall determine all the admissible weights $\lambda$ for  the rest of finite classical groups, that is, those such that
$\dim L(\lambda)<\ell^4$.  To this end, we will need to compute and bound a series of multiplicities $m_{\lambda}(\mu)$ for $\mu\preceq \lambda$.

\subsection{Type $B_\ell$}
We shall begin the desired classification of admissible weights by considering $G$ of  type $B_{\ell}$.
In order to compute the multiplicities involved in this section, we start by considering  $\{\alpha_1,\alpha_2,\cdots,\alpha_{\ell}\}$ the simple roots system (of type $B_\ell$) given by
\[
\alpha_1:=e_{\ell}\quad
\text{and} \quad\alpha_i:=e_{\ell-i+1}-e_{\ell-i+2},\quad\quad 2\leq i\leq \ell,
\]
where $\{e_1,e_2,\cdots, e_{\ell}  \}$ is an orthogonal basis for $E\cong \mathbb{R}^{\ell}$. 
We then note that  $(\alpha_1, \alpha_1)=1 $, $(\alpha_1, \alpha_2)= -1$, and for $i,j>1$ with $i\neq j$ we see that 
\begin{equation}
(\alpha_i, \alpha_j)= \begin{cases}
\phantom{-}2, & \text{ if }i=j\\
-1, & \text{ if }i=j\pm 1\\
\phantom{-}0, & \text{ otherwise }.\\
\end{cases}
\end{equation}

Then the  fundamental  weights $\lambda_i$ ($1 \leq i \leq \ell$), such that $(\lambda_i, \alpha_j) = \delta_{ij}$, are given by 
\begin{align*}
 \lambda_1  &= (1/2)(\alpha_{\ell}+2\alpha_{\ell-1}+\cdots+\ell\alpha_{1})\\
 \lambda_i &= \alpha_{\ell} + 2\alpha_{\ell-1}+ \cdots + (\ell-i)\alpha_{i+1}+(\ell-i+1)(\alpha_i+ \alpha_{i-1} + \cdots + \alpha_2+\alpha_1), \quad  1<i<\ell\\
 \lambda_{\ell} &= \alpha_\ell + \alpha_{\ell-1}+ \cdots + \alpha_2 + \alpha_1.
\end{align*}
In particular, we have that $(\lambda_i,\lambda_i) =
\ell-i+1.$ In addition, the half sum of the positive roots is 
\[
\delta=\frac{1}{2}\sum_{i=1}^{\ell}(\ell-i+1)(\ell+i-1)\alpha_{i}. 
\]
It is easy to verify 
\[
(\alpha_i,\delta) = \begin{cases}
\frac{1}{2}, & \text{ if }i=1\\
1, & \text{ otherwise }.\\
\end{cases}
\]

Let us start with the first example of admissibility for this type.
\begin{example}  Let $p$ be an odd prime and $V$ be the standard representation of $G(\overline{\F}_p)$ of dimension $2\ell + 1$.
For all $\ell \geq 12$, the weight $\lambda=4\lambda_{\ell}$ is  admissible. To see that, it suffices to show that $\dim \sym^{4}(V)$ is less or equal to $\ell^4$  since $L(4\lambda_\ell) \subseteq \sym^4(V)$. Indeed, we observe that 
\[
 \dim \sym^{4}(V)
 =\frac{(2\ell+1)(\ell+1)(2\ell+3)(\ell+2)}{6}
  =\ell^4\cdot \frac{(2+1/\ell)(1+1/\ell)(2+3/\ell)(1+2/\ell)}{6}\\
  \leq \ell^4,
\]
by verifying the inequality
 \[
(2+1/\ell)(1+1/\ell)(2+3/\ell)(1+2/\ell)\\
  \leq (2+1/12)(1+1/12)(2+3/12)(1+2/12)\leq 6,
 \]
 for any $\ell\geq 12$. Thus $\lambda=4\lambda_{\ell}$ is  an admissible weight.
\end{example}

Determining all the admissible weights $\lambda$ will imply to bound several multiplicities. We start with the following lemmas.
\begin{lemma}\label{multip_B_l2}
For all $\ell \geq 12$  and $p\neq 2$, a bound for the multiplicities of all dominant weights $\mu\preceq \lambda_{\ell-i}$, $i=2,3$, are given in the following table:

\begin{center}
\begin{tabular}{l| l r}
highest weight&dominant weight $\mu\preceq\lambda$ \quad &\quad  a bound of $m_{L(\lambda)}(\mu)$\\\hline
$\lambda_{\ell-2}$
& $\lambda_{\ell-2}$ & $1$ \\
& $\lambda_{\ell-1}$ & $1$ \\
& $\lambda_{\ell}$ & $\ell-1$ \\
&$ 0$&   $\ell$   \\\hline
$\lambda_{\ell-3}$
& $\lambda_{\ell-3}$ & $1$ \\
& $\lambda_{\ell-2}$ & $1$ \\
& $\lambda_{\ell-1}$ & $\ell-1$ \\
& $\lambda_{\ell}$ & $\ell$ \\
&$ 0$&   $\ell(\ell+1)/2$   \\\hline
\end{tabular}
\end{center}
\end{lemma}

\begin{proof}
We shall present the computations for $i=2$ and we  may remark that the case $i=3$ follows very similarly by resorting to the analogous arguments. So, we put $\lambda=\lambda_{\ell-2}$ and  $\mu=\lambda_{\ell-1}$.
Since 
$\mu=\lambda-(\alpha_{1}+\alpha_{2}+\cdots+\alpha_{\ell-2})$
we have to consider in the Freudenthal formula (\ref{FreudForm}), the positive roots of the form
\[
\alpha_i+\alpha_{i+1}+\cdots+\alpha_j, \qquad 1\leq i\leq j\leq \ell-2.
\]
Indeed, the formula for $m_{ \lambda}(\mu)$ involves the multiplicities of  $\mu+(\alpha_i+\cdots+\alpha_j)$ for any 
$ 1\leq i\leq j\leq \ell-2$.
 We then observe  
\begin{align*}
m_{ \lambda}(\mu)&=2\cdot \frac{\displaystyle{\sum_{1\leq i\leq j\leq \ell-2}m_\lambda(\mu+(\alpha_i+\cdots+\alpha_j))(\mu+(\alpha_i+\cdots+\alpha_j),\alpha_i+\cdots+\alpha_j)}}{( \lambda, \lambda)-(\mu,\mu)+2(\sum_{k=1}^{\ell-2}\alpha_k,\delta)}.\\
\end{align*}
First, by using $(\lambda_{\ell-2},\alpha_{k})=\delta_{(\ell-2)k}$, $(\alpha_1,\alpha_1)=1$ and $(\alpha_k,\alpha_{k+1})=-1$ we deduce   that 
\begin{align*}
( \lambda, \lambda)-(\mu,\mu)&=( \lambda_{\ell-2}, \lambda_{\ell-2})-\left(\lambda_{\ell-2}-\sum_{k=1}^{\ell-2}\alpha_k,\lambda_{\ell-2}-\sum_{k=1}^{\ell-2}\alpha_k\right)\\
&=2\left(\lambda_{\ell-2},\sum_{k=1}^{\ell-2}\alpha_k\right)-\left(\sum_{k=1}^{\ell-2}\alpha_k,\sum_{k=1}^{\ell-2}\alpha_k\right)\\
&=2\left(\lambda_{\ell-2},\alpha_{\ell-2}\right)-(\alpha_1,\alpha_1)-\sum_{k=2}^{\ell-2}(\alpha_k,\alpha_k)-2\sum_{k=1}^{\ell-3}(\alpha_k,\alpha_{k+1})\\
&=2-1-2(\ell-3)+2(\ell-3) =1.
\end{align*}
Thus 
\[
( \lambda, \lambda)-(\mu,\mu)+2(\sum_{k=1}^{\ell-2}\alpha_k,\delta)=1+2(1/2+(\ell-3))=2(\ell-2).
\]
Now, we claim that the positive roots $\alpha_1+\alpha_{2}+\cdots+\alpha_k $, with $1\leq k\leq \ell-2$, are the only ones that  contribute in the computation of $m_{\lambda}(\mu)$. Indeed, we can show that 
$m_{\lambda}(\mu+(\alpha_i+\cdots+\alpha_j))=0$,
for any $2\leq i\leq j\leq \ell-2$. To see that, we directly verify that $m_{\lambda}(\lambda-\alpha_{k})=0$ for any $1\leq k\leq\ell-3$ which implies inductively that $m_{\lambda}(\lambda-(\alpha_1+\cdots+\alpha_{i-1}))=0$, {\it i.e,} we prove the case $j=\ell-2$. Now, by assuming now that $j<\ell-2$  we have 
$\mu+(\alpha_i+\cdots+\alpha_j)=\lambda-(\alpha_1+\cdots+\alpha_{i-1})+
(\alpha_{j+1}+\cdots+\alpha_{\ell-2}).$
Then we can note that  the kind of  positive roots that contribute   in the Fruedental formula to compute the foregoing multiplicity are of the form $\alpha_{i'}+\cdots+\alpha_{j'}$ and $\alpha_{j}$, with $i\leq i'\leq j'\leq j$,  so  we  may deduce that 
$m_{\lambda}(\lambda-(\alpha_1+\cdots+\alpha_{i-1}+ \alpha_{j}))=0$ for any $i+1\leq j\leq \ell-2$.  By combining  these multiplicities  we can conclude that $m_\lambda(\lambda-(\alpha_1+\cdots+\alpha_{i-1})-
(\alpha_{j+1}+\cdots+\alpha_{t}))=0$, for $j+1\leq t\leq\ell-2$. In summary, we have proved that  
\begin{equation}\label{m_0}
m_\lambda(\mu+(\alpha_i+\cdots+\alpha_{j})) = 0, \quad 2\leq i\leq j<\ell-2.
\end{equation}
Next  by proving  
\begin{equation}\label{m_1}
m_\lambda(\mu+(\alpha_1+\cdots+\alpha_{k}))=1, \quad 1\leq k< \ell-2,
\end{equation} 
we can  deduce 
 \begin{align*}
 m_\lambda(\mu)&=2\cdot \frac{\displaystyle{\sum_{1\leq i\leq j\leq \ell-2}m_\lambda(\mu+(\alpha_i+\cdots+\alpha_j))(\mu+(\alpha_i+\cdots+\alpha_j),\alpha_i+\cdots+\alpha_j)}}{( \lambda, \lambda)-(\mu,\mu)+2(\sum_{k=1}^{\ell-2}\alpha_k,\delta)}\\
 &=2\cdot \frac{\displaystyle{\sum_{k=0}^{\ell-2}m_\lambda(\lambda_{\ell-2}-(\alpha_{k+1}+\cdots+\alpha_{\ell-2}))(\lambda-(\alpha_{k+1}+\cdots+\alpha_{\ell-2}),\alpha_1+\cdots+\alpha_k)}}{2(\ell-2)}\\
 &=2\cdot \frac{\ell-2}{2(\ell-2)} =1.
\end{align*}
Regarding  the proof of the  equation (\ref{m_1}) we proceed by induction.
We  begin by verifying 
\[
m_\lambda(\mu+(\alpha_1+\cdots+\alpha_{\ell-4}))=m_\lambda(\lambda-(\alpha_{\ell-3}+\alpha_{\ell-2}))=1
\]
thanks  to  
$m_\lambda(\lambda-\alpha_{\ell-2})=1$,     $m_\lambda(\lambda-\alpha_{\ell-3})=0$ and  $m_\lambda(\lambda)=1$. Now we assume by induction that
\[
m_\lambda(\mu+(\alpha_1+\cdots+\alpha_{k-1}))=m_\lambda(\lambda-(\alpha_{k}+\cdots +\alpha_{\ell-2}))=1,\quad 1\leq k<\ell-2,
\]
and we write   $\mu_{k}=\lambda-(\alpha_{k-1}+\alpha_k+\cdots +\alpha_{\ell-2})$.
By arguing as above and thanks to (\ref{m_0})  only the weights of the form  $\mu_{k}+(\alpha_{k-1}+\alpha_{k}+\cdots+\alpha_j)$ will contribute in the Fruedenthal formula of $m_\lambda(\lambda-(\alpha_{k-1}+\alpha_{k}+\cdots +\alpha_{\ell-2}))$. So, we  observe 
\[
(\mu_{k}+(\alpha_{k-1}+\cdots+\alpha_j),\alpha_{k-1}+\cdots+\alpha_j)=(\lambda-(\alpha_{j+1}+\cdots +\alpha_{\ell-2}),\alpha_{k-1}+\cdots+\alpha_j)=-(\alpha_{j+1},\alpha_j)=1.
\]
 Thus,
 \begin{align*}
 m_\lambda(\mu_j) &=2\cdot \frac{\displaystyle{\sum_{k=j}^{\ell-2}m_\lambda(\mu_{j}+(\alpha_{j}+\cdots+\alpha_k))(\mu_{j}+(\alpha_{j}+\cdots+\alpha_k),\alpha_j+\cdots+\alpha_k)}}{(\lambda,\lambda)-(\mu_j,\mu_j)-2(\sum_{k=j}^{\ell-2}\alpha_k,\delta)}\\
 &=2\cdot \frac{\displaystyle{\sum_{k=j}^{\ell-2}m_\lambda(\lambda-(\alpha_{k+1}+\cdots +\alpha_{\ell-2}))(\lambda-(\alpha_{k+1}+\cdots +\alpha_{\ell-2}),\alpha_j+\cdots+\alpha_k)}}{(\lambda,\lambda)-(\mu_j,\mu_j)-2(\sum_{k=j}^{\ell-2}\alpha_k,\delta)}\\
 &=2\cdot \frac{\ell-j-1}{2(\ell-j-1)} = 1.
\end{align*}
This completes the computation of $m_\lambda(\lambda_{\ell-1})=1$.

Next we shall compute $m_\lambda(\lambda_{\ell})$. By observing 
\[
\mu:=\lambda_{\ell}=\lambda_{\ell-2}-\sum_{i=1}^{\ell-2}\alpha_{i}-\sum_{i=1}^{\ell-1}\alpha_{i}
\]
we can infer that the only terms nonzero in the Freudenthal formula  depend on the following positive roots:
\begin{align*}
\alpha&=\sum_{i=1}^{\ell-2}\alpha_{i}+\sum_{i=1}^{\ell-1}\alpha_{i},\\
\alpha&=\alpha_{1}+\cdots+\alpha_k,\; \quad\:\;1\leq k\leq \ell-1,\\
\alpha&=\alpha_{k}+\cdots+\alpha_{\ell-1},\; \:\;2\leq k\leq \ell-1, \\
\alpha_{k,j}&=2(\alpha_{1}+\cdots+\alpha_k)+(\alpha_{k+1}+\cdots+\alpha_j),
\quad1\leq k\leq \ell-3\;\text{and } \ell-2\leq j\leq \ell-1.
\end{align*}
At this point, we need to remark that $m_{\lambda}(\mu+2(\alpha_{1}+\cdots+\alpha_{k}))(\mu+2(\alpha_{1}+\cdots+\alpha_{k}),(\alpha_{1}+\cdots+\alpha_{k}))$ also contributes in the Freudenthal formula of $m_{\lambda}(\mu)$, when $1\leq k\leq \ell-3$. According to the terms associated to the above positive roots  in  the foregoing formula,
we first compute that 
\[
m_{\lambda}\left(\mu+\sum_{i=1}^{\ell-2}\alpha_{i}+\sum_{i=1}^{\ell-1}\alpha_{i}\right)\cdot \left(\mu+\sum_{i=1}^{\ell-2}\alpha_{i}+\sum_{i=1}^{\ell-1}\alpha_{i},\sum_{i=1}^{\ell-2}\alpha_{i}+\sum_{i=1}^{\ell-1}\alpha_{i}\right)
=m_{\lambda}(\lambda)\left(\lambda,\sum_{i=1}^{\ell-2}\alpha_{i}+\sum_{i=1}^{\ell-1}\alpha_{i}\right)=1\cdot 2.
\]
Supposing that $\alpha_k=\alpha_{1}+\cdots+\alpha_k$, where $1\leq k\leq \ell-1$,  we directly verify that $(\mu+\alpha_k,\alpha_k)=1$. Now we shall prove that $m_{\lambda}(\mu+\alpha_k)=1$. Let us start by noticing that  the multiplicities for $k=\ell-i,\;i=1,2$ were computed in the above case, thus we have to compute $m_{\lambda}(\mu+\alpha_k)$ for  $k\leq \ell-3$. 
For the general case,  the Freudenthal formula involves computing many terms of the form $m_{\lambda}(\mu_k+\alpha))(\mu_k+\alpha,\alpha)$ for appropriate positive roots $\alpha$ we start by setting $\mu_k=\mu+\alpha_k$. 
Each of them implies in its turn to  compute many other multiplicities. Nevertheless, we introduce a case sufficiently  illustrative of these computations, say $m_{\lambda}(\mu_{\ell-3})$. Firstly, we see that 
\[
\mu_{\ell-3}=\lambda_{\ell-2}-(\alpha_{1}+\cdots+\alpha_{\ell-3}+2\alpha_{\ell-2}+\alpha_{\ell-1}).
\]
Immediately, we calculate
\[
(\lambda,\lambda)-(\mu_{\ell-3},\mu_{\ell-3})+2(\alpha_{\ell-1},\delta)=3-3+2\ell=2\ell.
\]
Now by appealing to similar arguments and strategies  as in the computing of $m_{\lambda_{\ell-2}}(\lambda_{\ell-1})$ we may deduce the only positive roots contributing to the formula to compute $m_{\lambda}(\mu_{\ell-3})$  are 
\[
\alpha_{1}+\cdots+\alpha_j,\quad\;\;1\leq j\leq \ell-1.
\]
Here we also argue  as in  $m_{\lambda}(\mu_{\ell-1})$  we 
may obtain 
\begin{align*}
m_{\lambda}(\mu+\alpha_1+\cdots+\alpha_{k})(\mu+\alpha_1+\cdots+\alpha_{\ell-1},\alpha_1+\cdots+\alpha_{k})&=1\cdot 1,\quad 1\leq k\leq \ell-1,\;\;k\neq \ell-3,\\
m_{\lambda}(\mu+\alpha_1+\cdots+\alpha_{\ell-3})(\mu+\alpha_1+\cdots+\alpha_{\ell-1},\alpha_1+\cdots+\alpha_{\ell-3})&=1\cdot 2. \\
\end{align*}
Therefore, we get
\[
m_{\lambda}(\mu)
=\frac{2[(\ell-2)\cdot 1+1\cdot 2]}{2\ell}=1.
\]
Regarding the term $m_{\lambda}(\mu+2\alpha_k)(\mu+2\alpha_k,\alpha_k)$, for $1\leq k\leq \ell-3$, we can proceed similarly as above to deduce 
each multiplicity $m_{\lambda}(\mu+2\alpha_k)$. At this point, routine computations allow to calculate the remaining components  to obtain  the multiplicities.

Finally, suppose that $\alpha_{k,j}=2(\alpha_{1}+\cdots+\alpha_k)+(\alpha_{k+1}+\cdots+\alpha_j)$, where
$1\leq k\leq \ell-3$ and  $\ell-2\leq j\leq \ell-1$.   We  compute $(\mu+\alpha_{k,j},\alpha_{k,j})$ as follows:
\begin{itemize}
\item For $j=\ell-1$ we have
\begin{align*}
(\mu+\alpha_{k,\ell-1},\alpha_{k,\ell-1})
&=(\lambda_{\ell-2}-(\alpha_{k+1}+\cdots+\alpha_{\ell-2}), 2(\alpha_{1}+\cdots+\alpha_k)+(\alpha_{k+1}+\cdots+\alpha_{\ell-1}))\\
&=(\lambda_{\ell-2},\alpha_{\ell-2})-2(\alpha_{k+1},\alpha_k)-(\alpha_{k+1},\alpha_{k+1})-(\alpha_{\ell-2},\alpha_{\ell-1}) = 2.
\end{align*}
Likewise we verify that $(\mu+\alpha_{k,\ell-2},\alpha_{k,\ell-2})=2$.
\item When $k\leq\ell-3$ we get 
\begin{align*}
(\mu+\alpha_{k,\ell-1},\alpha_{k,\ell-1})
&=(\lambda_{\ell-2}-(\alpha_{k+1}+\cdots+\alpha_{\ell-2}), 2(\alpha_{1}+\cdots+\alpha_k)+(\alpha_{k+1}+\cdots+\alpha_{\ell-1}))\\
&=(\lambda_{\ell-2},\alpha_{\ell-2})-2(\alpha_{k+1},\alpha_k)-(\alpha_{k+1},\alpha_{k+1})-(\alpha_{\ell-2},\alpha_{\ell-1})=2,
\end{align*}
and identically we get that $(\mu+\alpha_{k,\ell-2},\alpha_{k,\ell-2})=2$.
\end{itemize}

On the other hand, we may resort as above to determine the only positive roots $\alpha$ involving in the formula of $m_\lambda(\mu+\alpha_{k,j})$. Then, we recursively deduce that $m_{\lambda}( \lambda_{\ell-2}-(\alpha_{k+1}+\cdots+\alpha_{\ell-2}))=1$ as in the previous  case.
Since $1\leq k\leq \ell-3$ and  $\ell-2\leq j\leq \ell-1$, it follows that  the number of positive roots  $\alpha_{k,j}$ is 
\[
1\cdot 2+2\cdot 2+\cdots+2\cdot(\ell-3)=2\cdot (\ell-2)(\ell-3)/2= (\ell-2)(\ell-3).
\]
 In summary, the computations from the cases above provide respectively  the following summands to the formula of $m_{\lambda}(\mu)$:
\[
1\cdot 2,\quad (\ell-1)\cdot 1\cdot 1,\quad (\ell-3)\cdot 1\cdot 2,\quad (\ell-2)\cdot 1\cdot 2,\quad 1\cdot 2(\ell-2)(\ell-3)
\]
Then, we have that 
\[
m_{\lambda}(\mu)=\frac{2[2\cdot 1+(\ell-1)\cdot 1+2(\ell-2)+2(\ell-3)+2(\ell-2)(\ell-3)]}{2(2\ell-3)}=\ell-1.
\]
To obtain the multiplicity $m_\lambda(0)$, we observe that 
\[
\mu = 0 = \lambda - \left(\sum_{i=2}^{\ell-2} \alpha_{\ell-i} + \sum_{i=1}^{\ell-1} \alpha_{\ell-i}+\sum_{i=0}^{\ell} \alpha_{\ell-i}\right).
\]
This writing allows to apply the same strategy as above to determine all terms involved in the formula of $m_\lambda(0)$ and thereby to conclude the required multiplicity.  Then, the desired bound follows from (\ref{dimf}) and (\ref{dimf1})). 

The case $\lambda=4\lambda_{\ell}$ is  only longer, but completely similar than the ones considered above.  The natural first step is to compute the multiplicity of $\lambda_{\ell-1}+2\lambda_{\ell}$.
\end{proof}
\begin{lemma}\label{multip_B_l}
For all $\ell \geq 12$  and $p\neq 2$, a bound for the multiplicities of all dominant weights $\mu\preceq i\lambda_{\ell}$, $i=3,4$, are given in the following table:
\begin{center}
\begin{tabular}{l| l r}
highest weight&dominant weight $\mu\preceq\lambda$ & a bound of $m_{L(\lambda)}(\mu)$\\\hline
$3\lambda_{\ell}$
& $3\lambda_{\ell}$ & $1$ \\
& $\lambda_{\ell-1}+\lambda_{\ell}$  &   $1$ \\
&$\lambda_{\ell-2}$  &   $1$ \\
&$2\lambda_{\ell}$  &   $1$ \\
&$\lambda_{\ell-1}$  &   $1$ \\
&$\lambda_{\ell}$ &   $\ell$\\
&$ 0$&   $\ell$  \\
\hline
$4 \lambda_\ell$
& $4\lambda_\ell$ & $1$ \\
&  $\lambda_{\ell-1}+2\lambda_{\ell}$   &   $1$ \\
&  $2\lambda_{\ell-1}$   &   $1$ \\
&  $\lambda_{\ell-2}+\lambda_{\ell}$   &   $1$ \\
&  $\lambda_{\ell-3}$   &   $1$ \\
& $3\lambda_{\ell}$ & $1$ \\
& $\lambda_{\ell-1}+\lambda_{\ell}$  &   $1$ \\
&$\lambda_{\ell-2}$  &   $1$ \\
&$2\lambda_{\ell}$  &   $\ell$ \\
&$\lambda_{\ell-1}$  &   $\ell$ \\
&$\lambda_{\ell}$ &   $\ell$\\
&$ 0$&   $\ell(\ell-1)/2$  \\
\hline

\end{tabular}
\end{center}
\end{lemma}

\begin{proof} 

Let $\lambda=3\lambda_{\ell}$ and $\mu:=\lambda_{\ell-1}+\lambda_{\ell}=\lambda-\alpha_{\ell}$.
Then, we have that $(\lambda,\lambda)=9$ and  $(\mu,\mu)=5$.
In addition as $(\alpha_{\ell},\delta)=1$ and $(\lambda,\alpha_{\ell})=3(\lambda_{\ell},\alpha_{\ell})=1$, we have 
\[
m_{\lambda}(\mu)=2\cdot \frac{m_{\lambda}(\mu+\alpha_{\ell})(\mu+\alpha_{\ell},\alpha_{\ell})}{(\lambda,\lambda)-(\mu,\mu)+2(\alpha_{\ell},\delta)}=2\cdot \frac{1\cdot 3}{9-5+2\cdot 1}=1.
\]
Next,  compute the multiplicity of $\lambda_{\ell-2}=\lambda-(\alpha_{\ell-1}+2\alpha_{\ell})$, we observe that the positive roots involved in the  Fruedenthal formula are $\alpha_{\ell-1}$, $\alpha_{\ell}$, $2\alpha_{\ell}$ and $\alpha_{\ell-1}+\alpha_{\ell}$. Notice, at this point, that $\lambda$ is not included in the formula since  $\mu + \alpha_{\ell-1} + 2\alpha_\ell$ is not of the form $\mu + i\alpha$ for some positive root $\alpha$ by appealing to \cite[Plate II]{Bou02}. 

Set  $\mu:=\lambda-(\alpha_{\ell-1}+\alpha_{\ell})$.

\begin{enumerate}
 \item Set $\mu_{\ell-1}:=\mu+\alpha_{\ell-1}=\lambda-\alpha_{\ell}$.
 We first see that $m_{\lambda}(\mu_{\ell-1}+\alpha_{\ell})=1$ and $(\mu_{\ell}+\alpha_{\ell},\alpha_{\ell})=9$. In addition, we already  have $(\lambda,\lambda)=9$, and we verify that $(\mu_{\ell-1},\mu_{\ell-1})=5$. Thus, we conclude 
\[
m_{\lambda}(\mu_{\ell-1})=2\cdot \frac{m_{\lambda}(\mu_{\ell-1}+\alpha_{\ell})(\mu_{\ell-1}+\alpha_{\ell},\alpha_{\ell})}{(\lambda,\lambda)-(\mu_{\ell-1},\mu_{\ell-1})+2(\alpha_{\ell},\delta)}=1.
\]
 \item For  $\mu_{\ell}:=\mu+\alpha_{\ell}=\lambda-\alpha_{\ell-1}$,  we can conclude that $m_{\lambda}(\mu_{\ell})=0$ since $(\mu_{\ell}+\alpha_{\ell-1},\alpha_{\ell-1})=0$.
\end{enumerate}

Hence, by verifying that $(\mu+(\alpha_{\ell-1}+\alpha_{\ell}),\alpha_{\ell-1}+\alpha_{\ell})=3$, $(\mu,\mu)=5$ and $(\alpha_{\ell-1}+\alpha_{\ell},\delta)=2$
\begin{align*}
m_{\lambda}(\mu)&=2\cdot \frac{
m_{\lambda}(\mu_{\ell-1}+\alpha_{\ell-1})(\mu_{\ell-1}+\alpha_{\ell-1},\alpha_{\ell-1})+
m_{\lambda}(\mu_{\ell}+\alpha_{\ell})(\mu_{\ell}+\alpha_{\ell},\alpha_{\ell})
+m_{\lambda}(\lambda)(\lambda_{\ell},\alpha_{\ell-1}+\alpha_{\ell})}
{(\lambda,\lambda)-(\lambda_{\ell-2},\lambda_{\ell-2})+2(\alpha_{\ell-1}+\alpha_{\ell},\delta)} = 1.
\end{align*}
We continue with the computing of $m_{\lambda}(\lambda-2\alpha_{\ell})$. To this end we may deduce as above that $m_{\lambda}(\lambda-\alpha_{\ell})=1$ and next we shall compute the multiplicity $m_{\lambda}(\mu+2\alpha_{\ell})$.
We can verify that $(\mu+\alpha_{\ell}, \alpha_{\ell})=1$, $(\mu+2\alpha_{\ell},\alpha_{\ell})=3$ and $(\mu,\mu)=5$. Thus, we have
\begin{align*}
m_{\lambda}(\lambda-2\alpha_{\ell})&= 
 2\cdot \frac{m_{\lambda}(\lambda-\alpha_{\ell})(\lambda-\alpha_{\ell},\alpha_{\ell})+m_{\lambda}(\lambda)(\lambda,\alpha_{\ell})}{(\lambda,\lambda)-(\mu,\mu)+2(2\alpha_{\ell},\delta)} = 1.
\end{align*}
Going back to computing of $m_{\lambda}(\mu)$ with $\mu=\lambda_{\ell-2}=\lambda-(\alpha_{\ell-1}+2\alpha_{\ell})$ we have
\[
m_{\lambda}(\mu)=2\cdot \frac{(*)}{(\lambda,\lambda)-(\mu,\mu)+2(\alpha_{\ell-1}+2\alpha_{\ell},\delta)}
\]
where
\begin{align*}
    (*) & = (\mu + \alpha_{\ell},\alpha_\ell) m_\lambda(\mu+\alpha_\ell) + (\mu + \alpha_{\ell-1}) m_\lambda(\mu+\alpha_{\ell-1},\alpha_{\ell-1})  \\
    &+ (\mu + 2\alpha_{\ell},\alpha_\ell) m_\lambda(\mu+ 2\alpha_\ell) + (\mu + \alpha_{\ell-1} + \alpha_{\ell}, \alpha_{\ell-1}+\alpha_\ell) m_\lambda(\mu+\alpha_{\ell-1} + \alpha_\ell) \\
    & = (\mu + \alpha_{\ell}, \alpha_{\ell}) m_\lambda(\lambda-(\alpha_{\ell-1}+ \alpha_\ell)) + (\mu + \alpha_{\ell-1}, \alpha_{\ell-1}) m_\lambda(\lambda - 2\alpha_{\ell})  \\
    & +(\mu + 2\alpha_{\ell},\alpha_{\ell}) m_\lambda(\lambda -\alpha_{\ell-1}) + (\mu + \alpha_{\ell-1} + \alpha_{\ell}, \alpha_{\ell-1}+\alpha_\ell) m_\lambda(\lambda-\alpha_\ell). \\
    & = (\mu + \alpha_{\ell},\alpha_\ell) \cdot 1 + (\mu + \alpha_{\ell-1}, \alpha_{\ell-1}) \cdot 1 + (\mu + 2\alpha_{\ell},\alpha_\ell) \cdot 0 + (\mu + \alpha_{\ell-1} + \alpha_{\ell},\alpha_{\ell-1} + \alpha_{\ell}) \cdot 1 \\
     & = 2 + 2 + 0 + 2 = 6.
\end{align*}
Computing as above we see that 
\[
(\lambda,\lambda)-(\mu,\mu)+2(\alpha_{\ell-1}+2\alpha_{\ell},\delta)=9-3+6=12.
\]
Now, by using the multiplicities calculated above and computing the inner products, we have
\[
m_\lambda (\lambda_{\ell-2}) = \frac{2[1\cdot 2 +1\cdot 2+0\cdot 4+1\cdot 2]}{12} = 1. 
\]
The remaining  multiplicities $m_{\lambda}(\mu)$, for $\mu \preceq 2\lambda_{\ell}$, can be computed as in \cite{L}. 

Finally, the case $\lambda=4\lambda_{\ell}$ is longer but completely similar than the one above by calculating the multiplicity of $\lambda_{\ell-1}+2\lambda_{\ell}$ as first step.
\end{proof}

\begin{proposition}\label{DB} 
Let $G$ be as above  of type $B_\ell$ with  $\ell \geq 12$  and assume $p \neq 2$. The following weights 
\[
3\lambda_\ell\;4\lambda_\ell,\;\lambda_{\ell-2},\;\lambda_{\ell-3},\;\lambda_{\ell-1}+\lambda_\ell
\]
are  admissible weights. Additionally, $\lambda_1$ is also an admissible weight when $12 \leq \ell \leq 16$.
\end{proposition}

\begin{proof}
Suppose that $\ell\geq 12$  and $p \neq 2$.
The admissibility of the weights $3\lambda_\ell$   will be deduced by bounding  the dimension of 
$L(3\lambda_{\ell})$. The case of $4\lambda_\ell$ will follow similarly. To this end,  and according to  Eq. \ref{dimf1}, we have a bound  $m_{L(\lambda)}(\mu)$ for each dominant weight $\mu\preceq 3\lambda_{\ell}$ given in  Lemma \ref{multip_B_l}. On the other hand, observe that the stabilizers of the weights $\mu_i=i\lambda_{\ell}$, for $i=1,2,3$, are reflection groups of type $B_{\ell-1}$, which implies $l(\mu_i)=2\ell$. In addition,we also see that the stabilizers of $\lambda_{\ell-1}+\lambda_{\ell}$, $\lambda_{\ell-1}$ and $\lambda_{\ell-2}$ are reflection groups of type $B_{\ell-2}$, 
$B_{\ell-2}\times A_1$ and $B_{\ell-3}\times A_2$, so the length of its  orbits are  $4\ell(\ell-1)$,  $2\ell(\ell-1)$,  $4\ell(\ell-1)(\ell-2)/3$ respectively and clearly $l(0)=1$.
Then, by resorting to the multiplicities given in Lemma \ref{multip_B_l} and Eq. (\ref{dimf1}) we get
\begin{align*}
 \dim L(3\lambda_{\ell})
 &\leq  2\ell+4\ell(\ell-1)+
 4\ell(\ell-1)(\ell-2)/3+2\ell+2\ell(\ell-1)
 +2\ell^2+\ell\\ 
&=\ell(2\ell+1)(2\ell+5)/3,
\end{align*}
which implies  $3\lambda_{\ell}$ is admissible.
Similarly we can bound the dimensions of $V(\lambda_{\ell-2})$ and $V(\lambda_{\ell-3})$ by using Lemma \ref{multip_B_l2} and Eq. (\ref{dimf1}), and whence bounds for $L(\lambda_{\ell-2})$ and $L(\lambda_{\ell-3})$. Now, the weight $\lambda_{\ell-1}+\lambda_{\ell}$ is admissible by appealing to \cite[Table 3]{Ca}.
Finally, when  $12\leq \ell\leq 16$, we  infer that there not exist $\mu\preceq\lambda_1$ except $\lambda_1$ itself which implies that $\lambda_1$  is admissible since $l(\lambda_1) = 2^{\ell} \leq \ell^4$.
\end{proof}

In order to get a complete list of admissible weights of  type $B_\ell$, we will prove the following combinatorial lemma which will be very useful below. This result is perhaps well-known to the experts, but we give a proof for lack of a reference.

\begin{lemma}\label{comb}
The coefficients of 
\[
(2+x)^n = \sum_{k=0}^n \binom{n}{k}2^{n-k} x^k
\]
for $n \geq 12$ are  strictly increasing until $k= \lfloor \frac{n}{3} \rfloor$ and then they are  strictly decreasing.
\end{lemma}

\begin{proof}
Firstly, we set $t_k:=2^{n-k} \binom{n}{k}$ for each integer
$0\leq k\leq n$. Then
\[
t_k:=2^{n-k} \binom{n}{k} = 2 \cdot \frac{2^{n-(k+1)} n (n-1)\cdots(n-(k+1))}{k!}.
\]
We begin by observing  $t_0= 2^{n-1} \cdot 2 < 2^{n-1} \cdot n=t_1$ because $2<n$.  Notice also that 
\[
t_1=2^{n-1} n = 2^{n-2} n \cdot 2 < 2^{n-2} n \cdot \frac{n-1}{2}=t_2,
\]
thus $t_0<t_1<t_2$.
We may proceed as above to show $t_{k}<t_{k+1}$ when 
 $k < \lfloor \frac{n}{3} \rfloor$. It suffices to verify the limit case $k = \lfloor \frac{n}{3} \rfloor-1$.
In such a case,  from the inequality follows from 
\[
n-\left(\left\lfloor \frac{n}{3} \right\rfloor -1\right) \geq n- \frac{n}{3}  +1 = 2\frac{n}{3}+1\geq 2\left\lfloor\frac{n}{3}\right\rfloor,
\]
we can deduce that
\begin{align*}
t_{\lfloor \frac{n}{3} \rfloor -1} &=2^{n-\lfloor \frac{n}{3} \rfloor +1} \binom{n}{\lfloor \frac{n}{3} \rfloor -1}\\
&=
2^{n-\lfloor \frac{n}{3} \rfloor +1} \frac{n!}{(\lfloor \frac{n}{3} \rfloor -1)!(n-\lfloor \frac{n}{3} \rfloor )!(n-\lfloor \frac{n}{3} \rfloor +1)}\\
&=2^{n-\lfloor \frac{n}{3} \rfloor } \frac{n!}{(\lfloor \frac{n}{3} \rfloor -1)!(n-\lfloor \frac{n}{3} \rfloor )!}\cdot  \frac{2}{(n-\lfloor \frac{n}{3} \rfloor +1)}\\
&<2^{n-\lfloor \frac{n}{3} \rfloor } \frac{n!}{(\lfloor \frac{n}{3} \rfloor -1)!(n-\lfloor \frac{n}{3} \rfloor )!}\cdot  \frac{1}{\lfloor \frac{n}{3} \rfloor }\\
& =t_{\lfloor \frac{n}{3} \rfloor},
\end{align*}
thanks to the inequality $2/(n-\lfloor \frac{n}{3} \rfloor +1)<1/\lfloor\frac{n}{3}\rfloor$.
This proves that the coefficients of $(2+x)^n$ are strictly increasing until $k= \lfloor \frac{n}{3} \rfloor$.
On the other hand, by similar arguments, when $k > \lfloor \frac{n}{3} \rfloor$ we have that $2 \geq \frac{n - (k-1)}{k}$. Then, by observing,
\[
2^{n-(k+1)} \frac{n(n-1) \cdots (n-(k-2))}{(k-1)!} \cdot 2 \geq 2^{n-k} \frac{n(n-1) \cdots (n-(k-2))}{(k-1)!}\cdot \frac{(n-(k-1))}{k} 
\]
we show that the coefficients of $(2+x)^n$ are strictly decreasing when $k > \lfloor \frac{n}{3} \rfloor$. This finishes the proof.
\end{proof}
 Now we establish the main result concerning admissible weights of type $B_\ell$.
\begin{theorem}\label{PrinB}
Let $G$ be a finite group of Lie type $B_\ell$ over $\mathbb{\overline{F}}_p$ with $p>2$.
 If $\ell \geq 12$
 then the following list of weights
\[
0,\lambda_\ell,\;2\lambda_\ell,\; 3\lambda_\ell,\;4\lambda_\ell,\;\lambda_{\ell-1},\;\lambda_{\ell-2},\;\lambda_{\ell-3},\;\lambda_{\ell-1}+\lambda_\ell \mbox{ and } \lambda_1 \; (\mbox{when } 12 \leq \ell \leq 16)
\]
exhausts the admissible weights for $B_{\ell}$.
\end{theorem}
\begin{proof}
We shall proceed case-by-case  by arguing as in Theorem \ref{PrinA}.
Start   by setting 
\[
\lambda=a_1\lambda_1+a_2\lambda_2+\cdots+a_{\ell}\lambda_{\ell}.
\]
{\it Nullity of coefficients $-$ }
  $\lambda$ is admissible implies  $a_i=0$ for $1\leq i\leq \ell-4$.\\ 
Suppose  first that $a_i\neq 0$, for some  $ 2\leq i\leq \ell-4$. Then, the stabilizer of $\lambda$ is contained in a reflection group of type $B_{i-1}\times A_{\ell-i}$ - it is worth mentioning that  we are using the convention $B_1=A_1$ as in \cite{Con}. This  implies
\begin{align*}
\dim L(\lambda)\geq t_i:= \frac{2^{\ell}\ell !}{2^{i-1}(i-1)!(\ell-i+1)!}
&=2^{\ell-(i-1)}\binom{\ell}{i-1}.\\
\end{align*}
By numerical inspection we also  may verify for  each   $12\leq \ell\leq 15$ that
\begin{align*}
\dim L(\lambda)\geq
&2^{\ell-(i-1)}\binom{\ell}{i-1}>\ell^4,\quad\quad \text{ for }2\leq i\leq  \lfloor \frac{\ell}{3} \rfloor+1.
\end{align*} 
Now, for $\ell\geq 16$, we observe that 
\[
\dim L(\lambda) \geq t_i > t_0=2^{\ell}\geq \ell^4.
\]
by using   the $(i-1)$th-coefficients $t_i$ of $(2+x)^\ell$ are  strictly increasing for $2\leq i\leq  \lfloor \frac{\ell}{3} \rfloor+1$ and $2^{\ell}\geq \ell^4$.

On the other hand, we will demonstrate that $\dim L(\lambda)> \ell^4$ by using the sequence $(t_i)$ is decreasing for $i> \lfloor \frac{\ell}{3}\rfloor+1$. Indeed, we shall show that $t_{\ell-4} > \ell^4$ for $\ell\geq 12$ by observing
\begin{align*}
 t_{\ell-4} >\ell^4&\iff  2^5\frac{\ell(\ell-1)(\ell-2)(\ell-3)(\ell-4)}{5!} > \ell^4\\
  &\iff  2^5\frac{\ell^3(1-1/\ell)(1-2/\ell)(1-3/\ell)(\ell-4)}{2\cdot3\cdot 4\cdot 5} > \ell^3\\
&\iff  2^2\frac{(1-1/\ell)(1-2/\ell)(1-3/\ell)(\ell-4)}{3\cdot 5} > 1.\\
\end{align*}
Since $\ell\geq 12$ we can infer that $1-j/\ell \geq (12-j)/12$ for $1\leq j\leq 3$,
and $(\ell-4)/15\geq 8/15$.  By combining these facts we can infer 
\begin{align*}
2^2\frac{(1-1/\ell)(1-2/\ell)(1-3/\ell)(\ell-4)}{3\cdot 5}& =  2^2(1-1/\ell)(1-2/\ell)(1-3/\ell)(\ell-4)/15\\
&\geq 2^2\cdot \frac{11}{12}\cdot \frac{10}{12}\cdot \frac{9}{12}\cdot \frac{8}{15} > 1.
\end{align*}
With the foregoing  equivalence and inequality above, we can conclude that $t_{\ell-4}> \ell^4$ for  $\ell\geq  12 $.  This implies that 
$t_i\geq t_{\ell-4} > \ell^4$, for $i\geq  \lfloor \frac{\ell}{3} \rfloor+1$, when $\ell\geq 12$, by appealing to the sequence $(t_i)$ is decreasing in this rank.  We then conclude 
\[
\dim L(\lambda)\geq t_i\geq t_{\ell-4}> \ell^4,
\]
To deal the case $i=1$ we proceed as follows.  If $a_1\geq 2$ we deduce that $\lambda':=(\lambda-\alpha_{\ell-1}) \prec \lambda$ can be  written  as nonnegative linear combination of fundamental weights contains the term $\lambda_{2}$. Resorting as above   we have that 
\[
\dim L(\lambda)\geq\frac{\vert \mathcal{W}\vert }{ l(\lambda')} >\ell^4,
\]
{\it i.e.,} $\lambda$ is not admissible when $a_1\geq 2$.
By appealing to the same  arguments, $a_1=1$ implies $a_i=0$ for any $2\leq i\leq \ell$.
Hence we   reduce to 
\[
\lambda=\lambda_{1}\quad\quad\text{ or }\quad\quad
\lambda=a_{\ell-3}\lambda_{\ell-3}+a_{\ell-2}\lambda_{\ell-2}+a_{\ell-1}\lambda_{\ell-1}+a_{\ell}\lambda_{\ell}.
\]
By Proposition \ref{DB} we have that  $\lambda=\lambda_1$ is admissible when $12\leq \ell\leq 16$. On the other hand, since $\dim L(\lambda) = 2^{\ell} > \ell^4$ for $\ell > 16$ by \cite[Table 8]{Ca}, $\lambda_1$ is not admissible.      
Next we consider  $\lambda=a_{\ell-3}\lambda_{\ell-3}+a_{\ell-2}\lambda_{\ell-2}+a_{\ell-1}\lambda_{\ell-1}+a_{\ell}\lambda_{\ell}$. If $a_{\ell-3}\geq 2$, we produce a dominant weight $\lambda':=\lambda-\alpha_{\ell-3}$, whose writing as nonnegative linear combination of fundamental weights contains the term $\lambda_{\ell-4}$, but it cannot be possible due to the arguments  above. Thus, we reduce  the analysis to 
\[\lambda=\lambda_{\ell-3}+a_{\ell-2}\lambda_{\ell-2}+a_{\ell-1}\lambda_{\ell-1}+a_{\ell}\lambda_{\ell}\; \text{ or }\; \lambda=a_{\ell-2}\lambda_{\ell-2}+a_{\ell-1}\lambda_{\ell-1}+a_{\ell}\lambda_{\ell}.
\]
For the former case, {\it i.e.}, assuming $a_{\ell-3} = 1$,  we claim that $a_{\ell-j}=0$ for $j=0,1,2$.  In fact, if some of these coefficient were nonzero then the  stabilizer  $\rm{Stab}_\mathcal{W}(\lambda)$ of $\lambda$ were contained in a  reflection group of type   $B_{\ell-4}\times A_1$, $B_{\ell-4}\times A_1\times A_1$ or $B_{\ell-4}\times A_2$.
In any of those cases we can verify  that 
\[
\dim L(\lambda)\geq \frac{\vert \mathcal{W}\vert }{ l(\lambda')} >\ell^4.
\]
by arguing as cases above.

On the other hand,  when $\lambda=a_{\ell-2}\lambda_{\ell-2}+a_{\ell-1}\lambda_{\ell-1}+a_{\ell}\lambda_{\ell}$ we can prove that $a_{\ell-2}<2$. Otherwise, we observe that  $\lambda_{\ell-3}$ and $\lambda_{\ell-1}$ appear in the writing of  
$\lambda-\alpha_{\ell-2}$,  but this  case was  already covered previously. Thus,  we get  
\[
\lambda=\lambda_{\ell-2}+a_{\ell-1}\lambda_{\ell-1}+a_{\ell}\lambda_{\ell}\quad\text{ or }\quad 
\lambda=a_{\ell-1}\lambda_{\ell-1}+a_{\ell}\lambda_{\ell}.
\]
In both cases we may deduce that  $\lambda$ is not admissible if  $a_{\ell-i}>1$ with $i=0$ or $1$ by arguing as above. For instance, if $\lambda=\lambda_{\ell-2}+a_{\ell-1}\lambda_{\ell-1}+a_{\ell}\lambda_{\ell}$. with  $a_{\ell-1}\geq 2$ we can observe that the term $2\lambda_{\ell-2}$ appears in 
$\lambda':=\lambda-\alpha_{\ell-1}$, but  such a weight $\lambda'$ is not admissible by the previous case. It follows that  $\lambda=\lambda_{\ell-2}+a_{\ell}\lambda_{\ell}$ or  $\lambda=\lambda_{\ell-2}+\lambda_{\ell-1}+a_{\ell}\lambda_{\ell}$.  
In those cases we may proceed very similarly to rule out the weights when $a_{\ell}>1$. Thereby, we have to examine the admissibility of:
\[
\lambda_{\ell-2}+\lambda_{\ell-1},\; \lambda_{\ell-2}+\lambda_{\ell}\;\text{ and } \lambda_{\ell-2}+\lambda_{\ell-1}+\lambda_{\ell}
\]

For  $\lambda:=\lambda_{\ell-2}+\lambda_{\ell-1}$, by using  the identity,
\[
\alpha_{\ell-2}+\alpha_{\ell-1}=-\lambda_{\ell-3}+\lambda_{\ell-2}+\lambda_{\ell-1}-\lambda_{\ell}
\]
we  deduce that  $\mu:=\lambda_{\ell-3}+\lambda_{\ell} \prec \lambda_{\ell-2}+\lambda_{\ell-1}$. The same   calculations as above show that  
 the length of orbit of  $\mu$ is greater to $\ell^4$, and hence $\dim L(\lambda) > \ell^4$, {\it i.e.}, $\lambda$ is not admissible.
 On the other hand,  $\lambda:=\lambda_{\ell-2}+\lambda_{\ell}$ is not admissible by Example \ref{B21}.
While  for $\lambda:=\lambda_{\ell-2}+\lambda_{\ell-1}+\lambda_{\ell}$ one  notices that
$$
(\lambda_{\ell-2}+\lambda_{\ell-1}+\lambda_{\ell})-(\lambda_{\ell-3}+2\lambda_{\ell})=-\lambda_{\ell-3}+\lambda_{\ell-2}+\lambda_{\ell-1}-\lambda_{\ell}=
\alpha_{\ell-2}+\alpha_{\ell-1}
$$
which says  $\lambda_{\ell-3}+2\lambda_{\ell}\prec \lambda$. Since the stabilizer of $\lambda_{\ell-3}+2\lambda_{\ell}$ is of type $B_{\ell-4}\times B_2$, we have that 
 $$\dim L(\lambda)\geq l(\lambda_{\ell-3}+2\lambda_{\ell})
=2\ell(\ell-1)(\ell-2)(\ell-3) > \ell^4, $$
by observing $2(1-1/\ell)(1-2/\ell)(1-3/\ell)> 1$ for any $\ell\geq 12$. This implies that $\lambda=\lambda_{\ell-2}+\lambda_{\ell-1}+\lambda_{\ell}$ is not admissible either.

Finally, we consider the case $\lambda=a_{\ell}\lambda_{\ell}$. By appealing to Lemma \ref{DB}, the weights of the form $\lambda=a_{\ell}\lambda_{\ell}$ are admissible for  $a_{\ell}=1,2,3,4$. To complete this case we observe that none $\lambda$ with $a_{\ell}>4$ provides a new  admissible weight. To see that, it suffices to prove 
\[
\lambda_{\ell-4} \prec \lambda_{\ell-3}+\lambda_{\ell} \prec 5\lambda_{\ell} \preceq \lambda.
\]
forcing then $\dim L(\lambda)\geq  l(\lambda_{\ell-4}) > \ell^4$. Therefore, the weigths in Proposition \ref{DB} together with the weigths in \cite[Table 3]{Ca} exhausts the admissible weights which conclude the proof.
\end{proof}

Finally, we will deal with the rest of the classic groups of Lie type in the next subsection. 

\subsection{Type $C_\ell$ and $D_\ell$} 

\begin{theorem}\label{PrinCD}
 Let $G$ be a finite group of Lie type $C_\ell$ or $D_{\ell}$ over $\mathbb{\overline{F}}_p$ with $p>2$.
 If $\ell \geq 12$
 then the following list of weights
\[
0,\lambda_\ell,\;2\lambda_\ell,\; 3\lambda_\ell,\;4\lambda_\ell,\;\lambda_{\ell-1},\;\lambda_{\ell-2},\;\lambda_{\ell-3},\;\lambda_{\ell-1}+\lambda_\ell \mbox{ and }
\lambda_1, \lambda_2 \; (\mbox{when } 12 \leq \ell \leq 17 \mbox{ 
and } G \mbox{ of type }D_\ell ) \]
exhausts the admissible weights for $G$.
\end{theorem}

\begin{proof}
The tables in Appendix contain  all the information to determine bounds for the involved multiplicities   in the dimensions of the weights spaces $V(\lambda)$ and $L(\lambda)$. Hence, we may deduce that the aforementioned 
weights are effectively admissible. 
 Next we shall prove that the list is exhausted  respect to admissible ones. To this end, let
\[
\lambda=a_1\lambda_1+a_2\lambda_2+\cdots+a_{\ell}\lambda_{\ell}
\]
be a dominant weight and assume that $G$ is of type $C_\ell$. By arguing as in proof of Theorem \ref{PrinB}  (nullity of coefficients) we obtain that $\lambda$ admissible implies
\[
\lambda=\lambda_1\quad\quad\text{ or }\quad\quad\lambda=a_{\ell-3}
\lambda_{\ell-3}+a_{\ell-2}\lambda_{\ell-2}+a_{\ell-1}\lambda_{\ell-1}+a_{\ell}\lambda_{\ell}
\]
One here finds a little difference in determining the admissible weights with the case $B_{\ell}$, for $\ell=12,13,14,15$. Explicitly,  $\lambda=\lambda_1$ is  not admissible
 by noticing first $\lambda_3 \prec \lambda$ and that 
\[
 \dim L(\lambda)\geq l(\lambda_3) =2^{\ell-3}\ell(\ell-1)(\ell-2)>\ell^4.
\] 
Thus 
\[
 \lambda=a_{\ell-3}
\lambda_{\ell-3}+a_{\ell-2}\lambda_{\ell-2}+a_{\ell-1}\lambda_{\ell-1}+a_1\lambda_{\ell}.
\]
From this point, the rest of the route of $B_{\ell}$  works perfectly for the current case  $C_{\ell}$ by resorting to parallel arguments and  Appendix -Table \ref{table:1}.

Now,  assume that $G$ is of type $D_\ell$. By arguing as in proof of Theorem \ref{PrinB} too, we obtain that $\lambda$ admissible implies
\[
\lambda=\lambda_1,\quad\quad \lambda=\lambda_2 \quad\quad\text{or}\quad\quad\lambda=a_{\ell-3}
\lambda_{\ell-3}+a_{\ell-2}\lambda_{\ell-2}+a_{\ell-1}\lambda_{\ell-1}+a_{\ell}\lambda_{\ell}
\]
Contrary to the previous case, here, the weights $\lambda=\lambda_1$, $\lambda=\lambda_2$ are admissible for $12 \leq \ell \leq 17$, 
since the dimension of $L(\lambda_i)$ is $2^{\ell-1}$, for $i=1,2$,  by appealing to \cite[Lemma 2.3.2]{BGT}

Thereby, from
the dominant weights of the form $\lambda=a_{\ell-3}
\lambda_{\ell-3}+a_{\ell-2}\lambda_{\ell-2}+a_{\ell-1}\lambda_{\ell-1}+a_{\ell}\lambda_{\ell}$ on, the result  can be addressed as above by using  Appendix -Table \ref{table:2}.
\end{proof}


\section{Some dimension calculations and a potential application}\label{dimensiones}

It is well known that an explicit  formula for the  dimension of $L(\lambda)$ is unknown, or at least there is no known formula holding for an arbitrary $\lambda$. Nevertheless, Theorem \ref{PrinA} shows all these dimensions when $G$ is of type $A_{\ell}$ (See also \cite{Mar1}). On the other hand, computing such dimensions for certain weights spaces  seems to present major challenges when $G$ is not of type $A_{\ell}$,  where progress has been done case-by-case according to $G$ and very special weights $\lambda$ (See, for instance, \cite{Bur,BW,BGT,Ca,GS,Li,L,Mar2,Mc,Sei})  However, for $p$ sufficiently large, we are able to compute some of them  by  extending  arguments of \cite[\S 5]{L}, where dimension formulas for the admissible weights $0$, $\lambda_\ell$, $\lambda_{\ell-1}$ and $2 \lambda_\ell$ have been obtained.\\

\subsection{Type $C_\ell$ and $D_\ell$} 
 \begin{proposition}
Let $G$ of type $C_\ell$ over $\overline{\F}_p$ with $p>2$. Then,
\begin{enumerate}
 \item if $p \nmid (\ell-1)!$ then
    \[
  \quad \dim L(\lambda_{\ell-2})= \binom{2\ell}{3} - \binom{2\ell}{1}
    \]
   \item if $p \nmid \ell !$ then
 \[
 \quad\dim L(\lambda_{\ell-3})= \binom{2\ell}{4}-\binom{2\ell}{2},
 \]
    \item if $p \neq 3$ then
    \[
    \dim L(3 \lambda_{\ell}) = \binom{2\ell+2}{3} \qquad \mbox{ and } \qquad
    \dim L(4 \lambda_{\ell}) = \binom{2\ell+3}{4}.
    \]
\end{enumerate}
\end{proposition}

\begin{proof} For (1) and (2) we may argue similarly as in the proof of \cite[Theorem 5.1]{L}. First, for $\lambda_{\ell-2}$, we have the decomposition
\[
\Lambda^3(L(A_{\ell-1}))\oplus (\Lambda^2(L(A_{\ell-1}))\otimes L(A_{\ell-1})^*)\oplus ( L(A_{\ell-1})\otimes \Lambda^2(L(A_{\ell-1})^*)\oplus \Lambda^3(L(A_{\ell-1})^*)    
\]
However, in this case we can only ensure that $\Lambda^3(L(C_{\ell}))$ contains an irreducible constituent of dimension $\binom{2(\ell-1)}{3} - \binom{2(\ell-1)}{1}$, by applying induction, when $p \nmid ((\ell-1)-1)!$. Then by using the irreducibility of $\Lambda^2(V(A_{\ell-1}))\otimes V(A_{\ell-1})^*$ and $V(A_{\ell-1})\otimes \Lambda^2(V(A_{\ell-1})^*)$ for $p \nmid \ell-1$, we obtain our result.
Now, for $\lambda_{\ell-3}$, by using a similar decomposition as above, we may compute the dimension of $\Lambda^4(L(C_{\ell}))$, where the condition on $p$ comes from the induction on $\ell$ and the terms $\Lambda^2(L(A_{\ell-1}))\otimes \Lambda^2(L(A_{\ell-1})^*)$ and $\Lambda^3(L(A_{\ell-1}))\otimes L(A_{\ell-1})^*$.

For (3) we can also use a similar argument as above, however it follows more easily from \cite[Proposition 4.2.2.h]{Mc}.
\end{proof}

For type $D_\ell$, in addition to Lübeck's formulas \cite[Theorem 5.1]{L},  Burness, Ghandour, Testerman \cite[Lemma 2.3.2]{BGT}, Cavallin \cite[Table 9]{Ca} and Martínez \cite[Table 4]{Mar2} have obtained the following formulas
\[
\dim L(\lambda_1) = L(\lambda_2) = 2^{\ell-1}, \quad \dim L(3\lambda_{\ell})= \binom{2\ell+2}{3} - \binom{2\ell}{1} + \epsilon(\ell+1)(2\ell),
\]
and
\[
\dim L(\lambda_{\ell-i}) = \binom{2\ell}{i+1} \; \mbox{for} \; p \neq 2 \; \mbox{and} \; i=2,3.
\]
Moreover, if $p \neq 2,3$ and $p \nmid (\ell+2)!$, we can obtain the formula
 \[
 \quad\dim L(4 \lambda_{\ell})= \binom{2\ell +3}{4}-\binom{2\ell+1}{2}.
 \]
 as in the previous proposition.  The foregoing formulas, together with \cite[Theorem 5.1]{L}, allow us to compute the explicit dimension of $L(\lambda)$ (at least for sufficiently large primes $p$, explicitly $p>\ell+2$) for all admissible weights in Theorem \ref{PrinCD} except for $\lambda=\lambda_{\ell-1}+ \lambda_\ell$. However, it follows from \cite[\S 4]{Mar2} that
\begin{equation}\label{formulaCyD}
\dim L(\lambda) = \binom{2\ell+2}{3} \quad \mbox{for} \quad p \neq 2,3 \mbox{ and } p \nmid (2\ell-1)(2\ell+1)
\end{equation}
in both cases. Note that, this formula agrees with Lübeck's computations for small $\ell$ \cite[\S 6.32-6.48]{L}. 

For completeness, the dimension of $L(\lambda)$ for type $B_\ell$ may be found in \cite[\S 5]{L} for the weights $0,\; \lambda_{\ell},\; 2\lambda_{\ell},\; \lambda_{\ell-1}$; in \cite[Table 8]{Ca}  for $\lambda=\lambda_{\ell-a+1}$ with $a=3,4, \ell$; and in \cite{Mar2} for $3\lambda_{\ell}$ and $\lambda_{\ell-1}+\lambda_{\ell}$. For the remaining admissible weight $4 \lambda_{\ell}$ we believe that:
\begin{equation}\label{formulaB2}
  \dim L(4 \lambda_{\ell})=2\ell(2\ell+2)(2\ell+1)(2\ell+7)/24 - \epsilon(2\ell+5)\ell(2\ell+3) - \epsilon(2\ell+3)
  \end{equation}
The technique used in \cite{L} does not work in this case. Indeed, it seems to be appropriate to attempt a reverse approach, {\it i.e.,} embed our group into a group of a suitable type and use the information therein.

\subsection{Automorphic Galois representations}

Recall that a finite simple group of Lie type P$G$ in characteristic $p$ may be attached to a connected reductive simple algebraic group $G$ over $\overline{\F}_p$ of simply connected type and a Frobenius endomorphism $F:G(\overline{\F}_p) \rightarrow G(\overline{\F}_p)$ such that P$G \cong G^F/Z$, where $G^F  := \{ \gamma \in G(\overline{\F}_p) : F(\gamma) = \gamma \}$ is the group of fixed points of $F$ and $Z$ is the center of $G^F$. Then, the 
projective representations of P$G$ in characteristic $p$ are the same as linear representations of $G^F$ in characteristic $p$ (see \cite{Ste16}, p. 49, items (ix) and (x)), which can be constructed by restricting highest weight representations of $G$ to $G^F$. Consequently, knowing the representations of $G$, we will have available all the projective representations of $PG$.

 On the other hand, let $\rho_{p} : G_\Q \rightarrow \GO_n(\overline{\Q}_p) \subseteq \GL_n(\overline{\Q}_p)$ be a $p$-adic orthogonal Galois representation.
As we mentioned in the Introduction, examples of this kind of Galois representations can be obtained from the so-called regular algebraic, polarizable, cuspidal
automorphic representations $\pi$ of $\GL_n(\mathbb{A}_\Q)$, via the global Langlands correspondence \cite[\S 2.1]{BLGGT2014}.
In that case, the large image conjecture for automorphic Galois representations predicts that, for almost all primes $p$, the image of the projectivizations $\rho^{proj}_{p} : G_\Q \rightarrow {\rm PGO}_n(\overline{\F}_p)$ of the residual representations of $\rho_{p} : G_\Q \rightarrow \GO_n(\overline{\Q}_p)$ should contain $\POM^\pm_n(\F_{p^s})$ for some positive integer $s$. As we also pointed out at that moment, this conjecture is known in very few cases. However, if we impose certain conditions on the local components of $\pi$, we can at least bound the possible images of $\overline{\rho}_p$. For example, in \cite[Proposition 3.4]{Ze20}, the second author proved the following criterion concerning the image of certain families of orthogonal Galois representations.

\begin{proposition}\label{ref1}
Let $\mathcal{R} = \{ \rho_p \}_p$ be a family of $p$-adic orthogonal Galois representations $\rho_p : G_\Q \rightarrow \GO_n(\overline{\Q}_p)$ as in \cite[Theorem 2]{Ze20}. Then, for almost all primes $p$ (i.e. all but finitely many), the image of $\overline{\rho}_p^{proj}$ is contained in an almost simple group with socle a finite simple group of Lie type in characteristic $p$.
\end{proposition}

Then, a best knowledge   of the classification of irreducible representations of finite simple groups of Lie type can help us to improve Proposition \ref{ref1}. In particular, the main result of \cite{Ze20} (Theorem 1.1), states that for almost all primes $p$ and $n=4\varpi$ (where $\varpi$ is a prime number such that $17 \leq \varpi \leq 73$), the image of $\overline{\rho}_p^{proj}$ is equal to one of the following groups: 
\begin{equation}\label{realis}
\POM^\pm_n(\F_{p^s}), \; \PSO^\pm_n(\F_{p^s}), \; \PO_n^\pm(\F_{p^s}) , \; \PGO^\pm_n(\F_{p^s})
\end{equation} 
for some positive integer $s$. The idea of the proof is to look at the classification of projective irreducible representations of finite simple groups of Lie type P$G$ of dimension $n \leq 292$, given in \cite{L}, and prove that there are not irreducible representations of dimension $n$ other than those coming from a finite simple group of type $D_{2\varpi}$ for $p$ sufficiently large and $n = 4\varpi$ (see \cite[\S 2]{Ze20} and \cite{Ze19} for the precise conditions imposed on $p$ and $n$). This improving  had an interesting consequence concerning to the inverse Galois problem. 
Explicitly, \cite[Corollary 1.2]{Ze20} states that at least one of the groups in (\ref{realis}) can be realized as a Galois group over $\Q$ for infinitely many primes $p$ and infinitely many positive integers $s$.

In the same line of thought, using the strategy in Section 5 of \cite{Ze20} combined with Theorems \ref{PrinA}, \ref{PrinB} and \ref{PrinCD}, we may prove the following conditional improvement of Proposition \ref{ref1} for $n=4\varpi$ with $\varpi$ a prime number such that $79 \leq \varpi \leq 787$.

\begin{theorem}\label{princi}
Let $\varpi$ be a prime number such that $79 \leq \varpi \leq 787$, $n=4\varpi$, and $\mathcal{R} = \{ \rho_p \}_p$ be a compatible system of $n$-dimensional orthogonal Galois representations $\rho_p: G_\Q \rightarrow \GO_n(\overline{\Q}_p)$ as in \cite[Theorem 2]{Ze20}. Assuming that the formula (\ref{formulaB2}) is true, we have that for almost all $p$ the image of $\overline{\rho}_p^{proj}$ is equal to \begin{equation}\label{real1}
\POM^\pm_n(\F_{p^s}), \; \PSO^\pm_n(\F_{p^s}), \; \PO_n^\pm(\F_{p^s}) , \; \PGO^\pm_n(\F_{p^s}) 
\end{equation} 
or to a finite group $PG(\F_{p^s})$, where $G$ is of type $A_\ell$ ($\ell\leq 7$), $B_\ell$, $C_\ell$, $D_\ell$ ($\ell \leq 5$) or $G_2$, for some positive integer $s$.
\end{theorem}

The proof closely follows the ideas of \cite[Theorem 1.1]{Ze20}. Nevertheless, we will write a brief sketch of the proof for completeness. As we pointed out above,  projective representations of P$G$ in characteristic $p$ are the same as linear representations of $G^F$ in characteristic $p$, which can be constructed by restricting highest weight representations of $G$ to $G^F$. Thus, according to the discussion in \cite[\S 5]{Ze20} preceding the proof of Theorem 1.1 (which works for any dimension), it suffices to prove that the only irreducible $4\varpi$-dimensional representations $\overline{\rho}_p: G^F \rightarrow \Omega_{4\varpi}^\pm(\F_{q^\epsilon})$ that could occur, appear when $G$ is of type $A_\ell$ ($\ell \leq 7$), $B_\ell$, $C_\ell$, $D_\ell$ ($\ell \leq 5$), $G_2$ or $D_\ell$.

\begin{proof}
Assume that $G$ is of a Lie type different from $A_\ell$ ($\ell \leq 7$), $B_\ell$ , $C_\ell$, $D_\ell$ ($\ell \leq 5$) and $G_2$. In order to determine all irreducible representations in defining characteristic for groups of a fixed Lie type, we first need to compute all factorizations of $n$ into factors greater than 1. Such factorizations are $4\varpi = 4\cdot \varpi = 2 \cdot 2 \cdot \varpi = 2 \cdot 2 \varpi$. Appealing to Tables 6.6-6.53 of \cite{L} and the assumptions on $G$, there are no irreducible representations of dimension 2 or 4.  Thus, there are no tensor decomposable representations of dimension $4\varpi$. Moreover, according to Theorems \ref{PrinA}, \ref{PrinB} and \ref{PrinCD}; the dimension formulas in the previous section (including (\ref{formulaB2})); and Tables 6.6-6.53 of \cite{L},  the only irreducible tensor indecomposable self-dual representations of dimension $4\varpi$ are the natural representations of $C_{2\varpi}$ and $D_{2\varpi}$. However, the natural representation of $C_{2\varpi}$ is symplectic then the only irreducible $4\varpi$-dimensional orthogonal representation occurs when $G$ is of type $D_{2\varpi}$.
\end{proof}

Then, we will have the following interesting consequence concerning to the inverse Galois problem:

Once we prove that the formula (\ref{formulaB2}) is true for $p>2\ell+5$,  an interesting consequence concerning to the inverse Galois problem  can be  deduced as follows.
\begin{corollary}\label{Galois} 
Assuming that the formula (\ref{formulaB2}) is true for $p>2\ell+5$, we have that, for almost all prime $p$ and each $n$ of the form $4\varpi$ (where $\varpi$ is a prime number such that $79 \leq \varpi \leq 787$), at least one of the following groups: 
\begin{equation}\label{real2}
\PSL_{\ell+1}(\F_{p^s}), \; \PGL_{\ell+1}(\F_{p^s}), \; \PSU_{\ell+1}(\F_{p^s}) , \; \PGU_{\ell+1}(\F_{p^s}) \quad 2 \leq \ell \leq 7.
\end{equation}
\begin{equation}\label{real4}
G_2(\F_{p^s}), \; \POM_{2\ell+1}(\F_{p^s}), \; \PSO_{2\ell+1}(\F_{p^s}), \;  
\PSp_{2\ell}(\F_{p^s}), \; \PGSp_{2\ell}(\F_{p^s}), \quad 2 \leq \ell \leq 5.  
\end{equation} 
\begin{equation}\label{real1}
\POM^\pm_{2\ell}(\F_{p^s}), \; \PSO^\pm_{2\ell}(\F_{p^s}), \; \PO_{2\ell}^\pm(\F_{p^s}) , \; \PGO^\pm_{2\ell}(\F_{p^s}), \quad 2 \leq \ell \leq 5 \mbox{ or } \ell=2\varpi.  
\end{equation} 
 occurs as a Galois group over $\Q$ for infinitely many primes $p$ and infinitely many positive integers $s$. 
\end{corollary}

\begin{remark}
In the previous result, we include the possible occurrence of (\ref{real2}), (\ref{real4}), and (\ref{real1}), for $2 \leq \ell \leq 7$, as Galois groups of $\Q$. In the future, computational tools could be useful to compute all the representations of these groups for dimensions $n \leq 3148$ and obtain the exact analogue of \cite[Corollary 1.2]{Ze20}. It is worth pointing out  that such a result will be new for $\varpi > 181$ while for $\varpi \leq 181$ the foregoing result follows from \cite{L}.
\end{remark}

\begin{remark}
Finally we would like also to emphasize that if we assume the computations on the web page of Frack L\"ubek\footnote{http://www.math.rwth-aachen.de/~Frank.Luebeck/chev/WMSmall/index.html}  the Theorem \ref{princi} also works for $2477 \leq \varpi <10000$ and for $3037 \leq  \varpi <12155$ if we include groups of type $C_6$. In this case, the result will be new for $\varpi > 283$ and for $\varpi \leq 283$ follows from L\"ubek's computations.
\end{remark}

\begin{remark}
Besides this motivating application to Number Theory, there are many other interesting applications to Algebra itself. See, for example \cite{BL, HLM, LL, Le, Re, Ret}
\end{remark}

\section{Apendix}\label{A}
Next we  shall present the multiplicities and lengths of the dominants weights $\mu \preceq \lambda$, for each highest weight $\lambda$ listed in the table in Theorem \ref{PrinCD}. We may remark that  these tables provide all the information to bound the dimensions of  $V(\lambda)$ (and whence $L(\lambda)$) which implies admissibility of the foregoing weights $\lambda$.\\

In what follows we consider  orthogonal basis 
$\{e_1,e_2,\cdots, e_{\ell}  \}$ of $E\cong \mathbb{R}^{\ell}$. 

\subsection{Type $C_\ell$}
Start by considering  
$\{\alpha_1,\alpha_2,\cdots,\alpha_{\ell}\}$ the simple root system of type $C_\ell$ given by
$$
\alpha_1:= 2 e_{\ell}\quad
\text{and} \quad\alpha_i:=e_{\ell-i+1}-e_{\ell-i+2},\quad\quad 2\leq i\leq \ell,
$$
We then note that  $(\alpha_1, \alpha_1)=4$, $(\alpha_1, \alpha_2)= -2$ and  
\begin{equation}
(\alpha_i, \alpha_j)= \begin{cases}
\phantom{-}2, & \text{ if }i=j\\
-1, & \text{ if }i=j\pm 1\\
\phantom{-}0, & \text{ otherwise }
\end{cases}
\end{equation}
 for $i,j>1$. 
The  fundamental  weights $\lambda_i$ ($1 \leq i \leq \ell$), such that $(\lambda_i, \alpha_j) = \delta_{ij}$, are given by 
\begin{align*}
\lambda_1&=\alpha_{\ell} + 2\alpha_{\ell-1}+\cdots+(\ell-1)\alpha_2+\ell(1/2)\alpha_{1},\\
\lambda_i & =\alpha_{\ell}+2\alpha_{\ell-1}+\cdots+(\ell-i)\alpha_{i+1}+(\ell-i+1)(\alpha_{i}+\cdots+(1/2)\alpha_1),\quad  2<i<\ell,\\
\lambda_{\ell} & =\alpha_{\ell} + \alpha_{\ell-1}+\cdots\alpha_2+(1/2)\alpha_{1}.
\end{align*}
We note  that $(\lambda_i,\lambda_i) =\ell-i+1$. The half sum of the positive roots is 
\[
\delta=\sum_{i=1}^{\ell}\lambda_i=\ell\left(\frac{\ell+1}{2}\right)\cdot \frac{1}{2}\alpha_{1}+\sum_{i=2}^{\ell}(\ell-i+1)\left(\frac{\ell+i}{2}\right)\alpha_{i}.
\]
It is easy to verify that
\[
(\alpha_i,\delta) = \begin{cases}
2, & \text{ if }i=1\\
1, & \text{ otherwise}.\\
\end{cases}
\]
With the notation above in mind, the table below may be deduced by routing each step of the case $B_{\ell}$.
\begin{center}
\begin{table}[H]
\begin{tabular}{l| l c r}
highest weight & dominant weight $\mu \preceq \lambda$ &  \quad\quad bound of $m_{L(\lambda)}(\mu)$ & length of $\mathcal{W}$-orbit \\\hline
$3\lambda_{\ell}$
&$ 3\lambda_{\ell}$&   $1$ & $2\ell$ \\
&$\lambda_{\ell-1}+\lambda_{\ell}$&   $1$ & $4\ell(\ell-1)$ \\
& $\lambda_{\ell-2}$  &  $1$   & $2^2\ell(\ell-1)(\ell-2)/3$ \\
&$ \lambda_{\ell}$&  $\ell$  & $2\ell$  \\\hline
$4\lambda_{\ell}$
& $4\lambda_{\ell}$ & $1$ & $2\ell$ \\
& $\lambda_{\ell-1}+2\lambda_{\ell}$ & $1$  & $4\ell(\ell-1)$ \\
&$2\lambda_{\ell}$
&  $1$ & $2\ell$    \\
&$\lambda_{\ell-2}+\lambda_{\ell}$&  $1$  &$2^2\ell(\ell-1)(\ell-2)$ \\
&$\lambda_{\ell-3}$& $1$  &$2^4\ell(\ell-1)(\ell-2)(\ell-3)/24$ \\
& $2\lambda_{\ell}$ & $\ell$ & $2\ell$ \\
& $\lambda_{\ell-1}$ & $\ell$ & $2\ell(\ell-1)$ \\
& $0$ & $\ell(\ell+1)/2$ & $1$ \\
\hline
$\lambda_{\ell-2}$ & $\lambda_{\ell-2}$ & $1$ & $2^3\ell(\ell-1)(\ell-2)/6$ \\
& $\lambda_{\ell}$ &  $\ell-2$  & $2\ell$\\\hline
$\lambda_{\ell-3}$
& $\lambda_{\ell-3}$ & $1$ & $2^4\ell(\ell-1)(\ell-2)(\ell-3)/24$ \\
& $\lambda_{\ell-1}$ & $\ell-3$ & $2\ell(\ell-1)/6$ \\
&$0$&  $\ell(\ell-3)/2$  & $1$  \\
\hline
$\lambda_{\ell-1} + \lambda_{\ell}$
& $\lambda_{\ell-1} + \lambda_{\ell}$ & $1$ & $2^2\ell(\ell-1)$ \\
& $\lambda_{\ell-2}$ & $2$ & $2^2\ell(\ell-1)(\ell-2)/3$ \\
& $\lambda_{\ell}$ & $2(\ell-1)$ & $2\ell$ \\  
\hline
\end{tabular}
\caption{}
\label{table:1}
\end{table}
\end{center}

\subsection{Type $D_\ell$}

Start by considering  $\{\alpha_1,\alpha_2,\cdots,\alpha_{\ell}\}$ the simple root system of type $D_\ell$ given by
$$
\alpha_1:= e_{\ell-1} + e_{\ell}\quad
\text{and} \quad\alpha_i:=e_{\ell-i+1}-e_{\ell-i+2},\quad\quad 2\leq i\leq \ell.
$$
We then note that  $(\alpha_1, \alpha_1)=2$, $(\alpha_1, \alpha_2)= 0$, $(\alpha_1, \alpha_3)= -1$ and  
\begin{equation}
(\alpha_i, \alpha_j)= \begin{cases}
\phantom{-}2, & \text{ if }i=j\\
-1, & \text{ if }i=j\pm 1\\
\phantom{-}0, & \text{ otherwise}\\
\end{cases}
\end{equation}
for $i,j>1$.
The  fundamental  weights $\lambda_i$ ($1 \leq i \leq \ell$), such that $(\lambda_i, \alpha_j) = \delta_{ij}$, are given by 
\begin{align*}
 \lambda_1  &= (1/2)(\alpha_{\ell}+2\alpha_{\ell-1}+\cdots+(\ell-2) \alpha_{3} + (1/2)(\ell-2)\alpha_{2}+ (1/2)\ell  \alpha_{1}),\\
  \lambda_2  &= (1/2)(\alpha_{\ell}+2\alpha_{\ell-1}+\cdots+(\ell-2) \alpha_{3} + (1/2)\ell\alpha_{2}+ (1/2)(\ell-2)  \alpha_{1}),\\
  \lambda_i &= \alpha_{\ell} + 2\alpha_{\ell-1}+ \cdots + (\ell-i)\alpha_{i+1}+(\ell-i+1)(\alpha_i+ \alpha_{i-1} + \cdots + \alpha_3)+(1/2)(\ell-i+1)( \alpha_2+\alpha_1), \quad  2<i<\ell,\\
  \lambda_{\ell} &= \alpha_\ell + \alpha_{\ell-1}+ \cdots + \alpha_3 + (1/2)(\alpha_2 + \alpha_1).
\end{align*}
 
We observe that $(\lambda_i,\lambda_i) = \ell-i+1$ for $3 \leq i \leq \ell$ and $(\lambda_1,\lambda_1) = (\lambda_2,\lambda_2) = \ell/4$.
The half sum of the positive roots is 
\begin{align*}
\delta &= \frac{1}{4}\ell(\ell-1)(\alpha_1 + \alpha_2) + \frac{1}{2}\sum_{i=3}^{\ell}(\ell-i+1)(\ell+i-2)\alpha_{i}.
\end{align*}
It is easy to verify that
\[
(\alpha_i,\delta) = \begin{cases}
2 & \text{ if }i=1,2\\
1 & \text{ otherwise}.\\
\end{cases}
\]
With the notation above in mind, the table below may be deduced as in the previous subsection.

\begin{center}
\begin{table}[H]
\label{table:2}
\vspace{0.2 cm }
\begin{tabular}{l| l c r}
highest weight & dominant weight $\mu \preceq \lambda$ & \quad\quad a bound of $m_{L(\lambda)}(\mu)$ & length of $\mathcal{W}$-orbit \\\hline
$3\lambda_{\ell}$
& $3\lambda_{\ell}$  &   $1$ & $2\ell$ \\
&$ \lambda_{\ell-1}+\lambda_{\ell}$&   1 & $4\ell(\ell-1)$ \\
&$ \lambda_{\ell-2}$&   1 & $2^2\ell(\ell-1)(\ell-2)/3$ \\
&$ \lambda_{\ell}$&   $\ell-1$ & $2\ell$  \\\hline
$4\lambda_{\ell}$
& $4\lambda_{\ell}$ & $1$ & $2\ell$ \\
& $\lambda_{\ell-1}+2\lambda_{\ell}$ & $1$ & $4\ell(\ell-1)$ \\
& $2\lambda_{\ell-1}$ & $1$ & $2\ell(\ell-1)$ \\
& $\lambda_{\ell-2}+\lambda_{\ell}$ & $1$ & $2^2\ell(\ell-1)(\ell-2)$ \\
& $\lambda_{\ell-3}$ & $1$ & $2^2\ell(\ell-1)(\ell-2)(\ell-3)/24$ \\
& $2\lambda_{\ell}$ & $\ell$ & $2\ell$ \\
& $\lambda_{\ell-1}$ & $\ell$ & $2\ell(\ell-1)$ \\
&$ 0$&   $\ell(\ell-1)/2$ & 1  \\\hline
$\lambda_{\ell-2}$
& $\lambda_{\ell-2}$ & $1$ & $2^3\ell(\ell-1)(\ell-2)/6$ \\
& $\lambda_{\ell}$ & $\ell-1$ & $2\ell$\\\hline
$\lambda_{\ell-3}$
& $\lambda_{\ell-3}$ & $1$ & $2^4\ell(\ell-1)(\ell-2)(\ell-3)/24$ \\
& $\lambda_{\ell-1}$ & $\ell-2$ & $2\ell(\ell-1)$ \\
&$ 0$&   $\ell(\ell-1)/2$ & 1 \\\hline
$\lambda_{\ell-1} + \lambda_{\ell}$
& $\lambda_{\ell-1}+\lambda_{\ell}$ & $1$ & $2^2\ell(\ell-1)$ \\
& $\lambda_{\ell-2}$ & $2$ & $2^3\ell(\ell-1)(\ell-2)/6$ \\
& $\lambda_{\ell}$ & $2(\ell-1)$ & $2\ell$ \\\hline
\end{tabular}
\caption{}
\label{table:2}
\end{table}
\end{center}

\subsection*{Acknowledgments}
The second author was partially supported by the SECIHTI Grant CBF 2023-2024-224.

\bigskip
Luis Gutiérrez Frez \\
\textsc{Instituto de Ciencias Físicas y Matemáticas\\
Universidad Austral de Chile\\
Campus Isla Teja s/n, Valdivia\\
 Región de los Ríos, Chile}\\
 \emph{E-mail address:} \texttt{luis.gutierrez@uach.cl} \\

 \bigskip
Adrián Zenteno\\
\textsc{Universidad Autónoma Metropolitana - Unidad Iztapalapa\\
Departamento de Matemáticas, Edificio Anexo al T \\
Av San Rafael Atlixco No.186 Col.Vicentina \\ 
C.P. 09340, Iztapalapa, Mexico City, Mexico}\\
\emph{E-mail address: }\texttt{adrian.zenteno@xanum.uam.mx}

\end{document}